\documentclass[11pt]{article}
\usepackage[letterpaper,margin=1in]{geometry}
\usepackage[T1]{fontenc}
\usepackage[utf8]{inputenc}
\usepackage{lmodern,microtype}
\usepackage{amsmath,amssymb,amsthm,mathtools}
\usepackage{graphicx,booktabs,enumitem,xcolor}
\usepackage[hidelinks]{hyperref}
\usepackage{xspace}
\usepackage[capitalise,nameinlink,noabbrev]{cleveref}
\usepackage{comment}
\usepackage{tikz}
\usepackage{orcidlink}
\usepackage{pgfplots}
\pgfplotsset{compat=1.18}
\newtheorem{theorem}{Theorem}[section]
\newtheorem{proposition}[theorem]{Proposition}
\newtheorem{lemma}[theorem]{Lemma}
\newtheorem{corollary}[theorem]{Corollary}

\theoremstyle{definition}
\newtheorem{definition}[theorem]{Definition}

\newcommand{\Gnd}{\ensuremath{\mathcal G_{n,d}\xspace}}
\newcommand{\Sn}{\ensuremath{\mathfrak S_n}\xspace}
\newcommand{\Ov}{\operatorname{Ov}}
\newcommand{\GED}{\operatorname{GED}}
\newcommand{\Pp}{\mathbb P}
\newcommand{\Ee}{\mathbb E}
\newcommand{\dd}{\mathrm d}

\newcommand{\Prob}[1]{\mathbb P\left(#1\right)}
\newcommand{\E}[1]{\mathbb E \left[#1\right]}
\newcommand{\Var}[1]{\operatorname{Var}\left(#1\right)}

\newcommand{\eps}{\varepsilon}
\newif\ifdraftnotes
\draftnotestrue
\definecolor{draftblue}{RGB}{30,66,100}

\hypersetup{pdftitle={Maximal edge overlap of random regular graphs},
 pdfauthor={Andrea D'Ascenzo and Riccardo Michielan}}

\title{Maximum edge overlap between independent random regular graphs}
\author{Andrea D'Ascenzo\\Gran Sasso Science Institute, Italy\\
\texttt{andrea.dascenzo@gssi.it}
\orcidlink{0000-0001-5612-0798}
\and Riccardo Michielan\\Gran Sasso Science Institute, Italy\\
\texttt{riccardo.michielan@gssi.it}
\orcidlink{0000-0003-4642-8507}
}

\date{\today}
\begin{document}
\maketitle
\begin{abstract}
Optimising the vertex correspondence can substantially increase the number of edges shared by two independent realisations of the same random graph. We quantify this phenomenon for uniformly random $d$-regular graphs and derive corresponding results for their edge-based graph edit distance. For every fixed $d\geq 3$, we show that the maximum overlap is $n+\Theta_{\mathbb P}(n/\log n)$, where the leading term $n$ arises from Hamiltonicity, while the surplus captures the optimal order of the excess beyond it. We establish explicit upper and lower bounds on the surplus. The upper bound is obtained by counting common cores of minimum degree at least two, while the lower bound is established using spanning factors consisting of cycles of the same logarithmic length. For the special case $d=2$, we show that the maximum overlap is $n-(1+o_{\Pp}(1))\log n$. The additional overlap achievable is also studied by rotating aligned Hamilton cycles, using a contiguity argument. As a by-product, we determine the asymptotic maximum overlap between unions of independent random perfect matchings under cyclic shifts.
\end{abstract}

\section{Introduction}\label{sec:introduction}

Graph-structured representations offer a flexible way to model complex relationships and are therefore central to a broad spectrum of applications, ranging from semantic knowledge representation and molecular analysis to computer vision. In many of these settings, graphs rarely coincide exactly: they may differ in size, labelling, connectivity, or local structure while still encoding related information. This motivates the study of graph similarity, which replaces the yes-or-no notion of graph isomorphism with a graded measure of structural resemblance. By enabling meaningful comparisons between imperfect or partially aligned graphs, similarity measures provide a foundation for important learning tasks over graph data, including graph classification, clustering, pattern discovery, and anomaly detection.

Among the various approaches to quantify similarity, graph overlap is a widely used and intuitive measure. Rather than requiring an exact structural correspondence, graph overlap measures the extent to which the two graphs share common structural elements, such as vertices, edges, or substructures. This makes graph overlap particularly suitable for applications in which graphs may be only partially similar or may contain noise, missing information, or additional structures.

We quantify the similarity between two graphs by the number of shared edges they have. More precisely, we consider two graphs defined on the same vertex set. 

\begin{definition}[Graph overlap]
Let \(G_1\) and \(G_2\) be two graphs with a common vertex set \(V\), where \(n=|V|\). For a permutation \(\pi\in\Sn\) of \([n]=\{1,\ldots,n\}\), the graph overlap of \(G_1\) and \(G_2\) under \(\pi\) is defined by
\[
    \Ov_{\pi}(G_1,G_2)
    =
    \bigl|\pi(E(G_1))\cap E(G_2)\bigr|.
\]
The maximum overlap is then given by
\begin{equation}\label{eq:objective}
    M_n
    =
    \max_{\pi\in\Sn}
    \Ov_{\pi}(G_1,G_2).
\end{equation}
\end{definition}

The quantity \(M_n\) is closely related to the graph edit distance (GED)~\cite{DBLP:journals/tsmc/SanfeliuF83}. In general, the GED between two graphs, possibly equipped with vertex and edge labels, is defined as the minimum cost of an edit sequence that transforms one graph into a graph isomorphic to the other. The allowed edit operations may include inserting, deleting, or substituting vertices and edges, with each operation assigned a non-negative cost \cite{gao2010survey}.

For unlabelled graphs on the same vertex set, when only edge insertions and deletions are permitted, and every operation has unit cost, the edge-based GED is given by
\begin{equation}\label{eq:ged-overlap}
    \GED_{E}(G_1,G_2)
    =
    \min_{\pi\in\Sn}
    \bigl|\pi(E(G_1))\triangle E(G_2)\bigr| = |E(G_1)|+|E(G_2)|-2M_n.
\end{equation}
In this paper, we assume that $G_1$ and $G_2$ are independent copies of the same graph distribution, so that $M_n$ represents a notion of self-similarity within the same graph ensemble. Both graph overlap and edit distance have been investigated in the setting of independent Erd\H{o}s--R\'enyi random graphs. 

\paragraph{Related work}

Comparisons of independent realisations of the same random graph model have been investigated from several perspectives. Hartle et al.~\cite{HartleEtAl2020} introduced the notion of \emph{within-ensemble graph distance}, studying the typical distance between two independent samples from a common graph ensemble. Their analysis covers several random graph models, including Erd\H{o}s--R\'enyi graphs, random geometric graphs, Watts--Strogatz networks, configuration models, and preferential attachment networks, using a variety of graph distance measures. This provides a general framework for understanding how structural variability within an ensemble is reflected by different notions of graph similarity.

From a mathematical perspective, Bollobás and Riordan~\cite{Bollobás_Riordan_2009,DBLP:journals/rsa/BollobasR11} investigated metrics on sparse graphs, including edit and cut distances, with particular attention to random graphs having a linear number of edges. In the Erd\H{o}s--R\'enyi model $\mathcal G(n,p)$ with $p=c/n$, their results show that, for fixed $c>1$, two independent realisations have an edit distance bounded below by a positive constant times $n$ with high probability. Their framework also addresses more general inhomogeneous sparse random graph models. More recently, Ding, Du, and Gong~\cite{DingDuGong2024} studied the maximum overlap between two independent Erd\H{o}s--R\'enyi graphs in denser regimes. For $p=n^{-\alpha}$ with $\alpha\in(1/2,1)$, they established that
\[
    M_n \sim \frac{n}{2\alpha-1}
\]
with high probability, and provided a polynomial-time approximation scheme achieving an arbitrarily close multiplicative approximation to the optimum. This shows that optimising the vertex correspondence can produce a linear overlap even when the two graphs are independently generated.

The algorithmic aspects of this problem were further investigated by Du, Gong, and Huang~\cite{DuGongHuang2025}, who identified a phase transition around
\[
    p_c=\sqrt{\frac{\log n}{n}}.
\]

They established polynomial-time approximation schemes in the sparse regime and identified algorithmic limitations in the dense regime, including hardness results for online and greedy local iterative algorithms. Their overlap analysis also provides evidence of a broader statistical-computational gap, although the hardness of arbitrary polynomial-time algorithms remains open.

In this paper, we investigate the maximum overlap within the ensemble of uniform random $d$-regular graphs. Unlike Erd\H{o}s--R\'enyi graphs, these graphs satisfy a hard degree constraint, which introduces dependencies between their edges. Moreover, for fixed $d\geq3$, the presence of Hamilton cycles provides a natural linear baseline for the overlap. Our aim is to determine not only this leading-order behaviour, but also the order of the additional overlap achievable by optimising over vertex permutations.

\subsection{Labelled versus unlabelled overlap}
We focus on the set $\Gnd$ of simple $d$-regular graphs on $[n]$, where $d$ is fixed and $nd$ is even. We denote by $\mathcal U_{n,d}$ the uniform probability measure on $\Gnd$. Throughout the paper, $G_1$ and $G_2$ are independent samples from $\mathcal U_{n,d}$, and all asymptotic statements are understood along admissible values of $n$.

To illustrate the effect of optimising over vertex permutations, we first consider the overlap under a fixed labelling. Let
\[
    X_n := |E(G_1)\cap E(G_2)|
\]
be the number of common edges. Since each unordered pair of vertices is an edge of a uniform $d$-regular graph with probability $d/(n-1)$, independence gives
\[
    \E{X_n} = \binom{n}{2}\left(\frac{d}{n-1}\right)^2 = \frac{nd^2}{2(n-1)} = \frac{d^2}{2}+o(1).
\]

In fact, using standard estimates for random regular graphs (see, e.g., \cite[Chapter~9]{janson2000randomgraphs}), a factorial-moment argument shows that, for every fixed $d$, $X_n$ converges in distribution to a Poisson random variable with mean $d^2/2$.
Thus, the overlap under a fixed labelling remains of constant order as $n$ grows, and in particular vanishes with a probability bounded away from zero.

Allowing vertex relabelings changes the picture. For any two graphs $G_1$ and $G_2$ on the same vertex set, consider a uniformly random permutation $\pi\in\Sn$. Since the image of each fixed edge is uniformly distributed over the $\binom{n}{2}$ unordered vertex pairs, we have
\[
    \mathbb{E}_{\pi}\bigl[\Ov_{\pi}(G_1,G_2)\bigr] = \frac{|E(G_1)|\,|E(G_2)|}{\binom{n}{2}}.
\]
Consequently, at least one permutation achieves an overlap no smaller than this average. Since every $d$-regular graph has exactly $nd/2$ edges, we obtain the deterministic bound
\[
    M_n \geq \left\lceil\frac{nd^2}{2(n-1)}\right\rceil.
\]

Although this elementary argument only yields a constant lower bound, optimising over vertex permutations has a much stronger effect: as we show below, the maximum overlap between two independent random $d$-regular graphs is typically of linear order in $n$.

\subsection{Main results}
Our objective is to characterise the maximum overlap between two independent uniform random $d$-regular graphs, for each fixed $d\geq1$. We will show that the behaviour changes significantly between the regimes $d=1$, $d=2$, and $d\geq3$ as summarised in the following table:

\begin{table}[htbp]
\centering
\label{tab:overlap-summary}
\begin{tabular}{|c|c|c|}
\hline
$\quad$ Degree $\quad$ & Maximum overlap $M_n$ & $\GED_E(G_1,G_2)$ \\
\hline
$d=1$ & $n/2$ & $0$ \\
$d=2$ & $n-(1+o_{\Pp}(1))\log n$ & $(2+o_{\Pp}(1))\log n$ \\
$d\geq3$ & $n+\Theta_{\mathbb P}(n/\log n)$ & $(d-2)n-\Theta_{\mathbb P}(n/\log n)$ \\
\hline
\end{tabular}
\end{table}

When $d=1$, both graphs are perfect matchings, and any two perfect matchings on the same vertex set are isomorphic. Consequently, all $n/2$ edges can be aligned, yielding $M_n=n/2$ and $\GED_E(G_1,G_2)=0$ deterministically.

The case $d=2$ is more interesting. Here, each graph is a disjoint union of cycles, and although two independent realisations are typically not isomorphic, their cycle decompositions can be aligned so that only a logarithmic number of edges remain unmatched.

\begin{theorem}[Overlap in degree two]\label{thm:overlap-d2}
Let $d=2$. Then
\[
\frac{n-M_n}{\log n}\xrightarrow{\Pp}1.
\]
More precisely, with $c(G)$ denoting the number of cycles of $G$,
\[
M_n=n-c(G_1)-c(G_2)+O_{\Pp}(1).
\]
\end{theorem}
Since both graphs have exactly $n$ edges, this also implies that $\GED_E(G_1,G_2)=(2+o_{\Pp}(1))\log n$. Thus, unlike the case $d=1$, two independent random $2$-regular graphs are typically not isomorphic, but they can be made almost identical by an appropriate relabelling, with only a logarithmic number of edge modifications.

For every fixed $d\geq3$, the behaviour is fundamentally different. A classical result of Robinson and Wormald~\cite{RobinsonWormald1994} guarantees that a uniform random $d$-regular graph is Hamiltonian with high probability. Consequently, two independent realisations contain Hamilton cycles that can be perfectly aligned through a suitable vertex permutation, yielding $M_n\geq n$ with high probability. On the other hand, a first-moment argument shows that the maximum overlap cannot exceed $(1+\varepsilon)n$ with high probability for any fixed $\varepsilon>0$.

\begin{theorem}[First-order overlap]\label{thm:first-order-overlap}
For every fixed $d\geq3$ and every $\varepsilon>0$,
\[
    \lim_{n\to\infty}\mathbb P\left(\left|\frac{M_n}{n}-1\right|\leq\varepsilon\right)=1.
\]
\end{theorem}

This determines the leading-order behaviour of the maximum overlap. The next question is whether an appropriate relabelling can systematically improve upon the overlap of $n$ edges provided by aligning two Hamilton cycles. To quantify this improvement, we introduce the \emph{surplus overlap}
\[
    S_n:=M_n-n.
\]
Our main results identify the correct asymptotic order of this surplus and give explicit degree-dependent constants.  Set
\[
 a_d:=\log(d-1)+\frac{d-2}{2}\log\left(\frac{d-2}{d}\right),
 \qquad
 b_d:=d\log(d-1)-(d-2)\log(d+1).
\]
We first establish an upper bound showing that, for any fixed $d\geq3$, the surplus cannot exceed order $n/\log n$ with high probability.

\begin{theorem}[Second-order upper bound]\label{thm:second-order-upper}
For every fixed $d\geq3$ and every $\varepsilon>0$,
\[
    \lim_{n\to\infty}\mathbb P\left(S_n\leq (b_d+\varepsilon)\frac{n}{\log n}\right)=1.
\]
\end{theorem}

The proof again uses a first-moment argument, but now after a structural reduction.  Any overlap of size $n+s$ contains a common subgraph (a \textit{core}) of minimum degree at least two and graph excess exactly $s$.  Counting these common cores and optimising over their number of vertices gives the constant $b_d$; the first moment is maximised by cores occupying asymptotically $d/(d+1)$ of the vertices.

For the lower bound, we first fix a Hamilton cycle in each graph and consider a permutation that perfectly aligns the two cycles, yielding an overlap of at least $n$ edges. Cyclically rotating this permutation preserves the alignment of all $n$ Hamiltonian edges, while changing the correspondence between the remaining edges. We exploit these $n$ possible rotations to search for additional overlaps beyond the Hamiltonian baseline. Using contiguity arguments, the residual edges can be analysed through a suitable random matching model. We then identify local configurations witnessing additional overlaps and control the dependencies among different rotations by selecting witnesses with distinct cyclic lengths. This yields a surplus of order $\log n/\log\log n$ with high probability.

To reach the larger scale $n/\log n$, we need to pass from a single Hamilton cycle to a spanning collection of cycles of the same logarithmic length $k$. 
\begin{theorem}[Second-order lower bound]\label{thm:second-order-lower}
For every fixed $d\geq3$ and every $\varepsilon>0$,
\[
    \lim_{n\to\infty}\mathbb P\left(S_n\geq (a_d-\varepsilon)\frac{n}{\log n}\right)=1.
\]
\end{theorem}

In this case, we use a different probabilistic mechanism. At the correct length scale $k$, the expected number of realised $k$-cycle factors is sub-exponential in $n$, which allows us to prove two properties simultaneously for every such factor. We can then choose the spanning factor after examining the graph and, consequently, align groups of cycles as stars. This procedure gains one additional common edge per factor cycle, asymptotically, and at the same time, no cycle-factor contiguity argument is needed.

\begin{figure}[t]
    \centering
    \begin{tikzpicture}
    \begin{axis}[
        width=0.5\textwidth,
        height=0.3\textwidth,
        xlabel={$d$},
        ylabel={second-order constant},
        xmin=2.5,
        xmax=20.5,
        ymin=0,
        xtick={3,5,10,15,20},
        tick align=outside,
        axis line style={black!75},
        tick style={black!75},
        grid=major,
        grid style={black!10},
        legend style={
            draw=none,
            fill=none,
            at={(0.03,0.97)},
            anchor=north west,
            cells={anchor=west}
        },
        legend cell align=left,
    ]

    \addplot[
        very thick,
        solid,
        black,
        mark=*,
        mark size=1.65pt,
        mark options={solid,fill=white}
    ]
    expression[domain=3:20,samples=18]
    {ln(x-1) + (x-2)/2*ln((x-2)/x)};
    \addlegendentry{$a_d$}

    \addplot[
        very thick,
        dashed,
        black,
        mark=square*,
        mark size=1.55pt,
        mark options={solid,fill=black}
    ]
    expression[domain=3:20,samples=18]
    {x*ln(x-1) - (x-2)*ln(x+1)};
    \addlegendentry{$b_d$}

    \end{axis}
    \end{tikzpicture}
    \caption{Lower- and upper-bound constants $a_d$ and $b_d$ in the second-order overlap bounds in Theorems \ref{thm:second-order-upper}-\ref{thm:second-order-lower}, for two independent random $d$-regular graphs.}
    \label{fig:second-order-constants}
\end{figure}

Together, the bounds establish the correct asymptotic order of the surplus, $S_n=\Theta_{\mathbb P}(n/\log n)$. Figure \ref{fig:second-order-constants} depicts the distance between the upper and lower constants. For $d=3$, they are
\[
 a_3=\frac12\log\frac43\approx0.14384,
 \qquad
 b_3=\log2\approx0.69315.
\]
More generally, $a_d=\log d-1+O(d^{-2})$ while $b_d=2\log d-2+2/d+O(d^{-2})$, so the ratio of the coefficients tends to $2$, as $d$ grows.

\subsection{Further results}

As a by-product, our analysis also yields an independent result concerning the maximum overlap between unions of random perfect matchings under cyclic shifts.

\begin{theorem}[Overlap of independent perfect matchings]\label{thm:cyclic-matching-unions}
Fix positive integers $q_1,q_2$ and let $n$ be even. For $i\in\{1,2\}$, let $\mathcal M_{i,1},\ldots, \mathcal M_{i,q_i}$ be mutually independent uniform perfect matchings on $\mathbb Z_n$, with all matchings independent across the two families, and define
\[
    R_i=\bigcup_{j=1}^{q_i}\mathcal M_{i,j}.
\]
Then, for every $\varepsilon>0$,
\[
    \lim_{n\to\infty}\mathbb P\left(
    \left|
    \frac{\displaystyle\max_{s\in\mathbb Z_n}|\rho_s(R_1)\cap R_2|}
    {\log n/\log\log n}
    -1
    \right|\leq\varepsilon
    \right)=1,
\]
where $\rho_s(x)=x+s\pmod n$ denotes the cyclic shift by $s$.
\end{theorem}

This result identifies the precise asymptotic scale of the largest overlap achievable through cyclic shifts of unions of independent random matching. It also provides the probabilistic ingredient underlying our analysis of rotations of aligned Hamilton cycles. Its proof is deferred to \cref{app:hamiltonian-technical} in the appendix.

A further ingredient in our second-order lower bound is the existence of $k$-cycle factors in random regular graphs, where $k$ is of order $\log n$. Kahn, Lubetzky, and Wormald~\cite{KahnLubetzkyWormald2017} established the sharp threshold in random \textit{cubic} graphs and explicitly state that their proof extends, with very few modifications, to every fixed degree $d\ge3$. In \cref{sec:preliminaries}, we extract the degree-dependent threshold constant from the first moment.

\subsection{Organisation of the paper}
\cref{sec:preliminaries} collects model definition, Hamilton-cycle contiguity, and the cycle-factor threshold used later.  \cref{sec:easy_section} contains the proofs for $d=1,2$ and first-order overlap estimates for fixed degree $d \geq 3$.  \cref{sec:upper} proves the second-order upper bound by reducing a large overlap to a common minimum-degree-two core.  \cref{sec:hamiltonian-alignment} studies the weaker but more explicit surplus obtainable by rotations of aligned Hamilton cycles, using a planted-model contiguity argument. \cref{sec:lower} proves the matching-order lower bound from logarithmic cycle factors by an argument working uniformly over all cycle factors and a deterministic star construction.  We conclude with a discussion of the remaining constant gap and related open questions.

\section{Preliminaries}\label{sec:preliminaries}

We use the notation $\Gnd$ and $\mathcal U_{n,d}$ introduced in \cref{sec:introduction}. Unless otherwise specified, $d$ is fixed, and all asymptotic statements are understood along admissible values of $n$, namely those for which $nd$ is even. For two independent samples $G_1,G_2\sim\mathcal U_{n,d}$, recall that
\[
M_n=\max_{\pi\in\Sn}|\pi(E(G_1))\cap E(G_2)|.
\]
Note that, in what follows, all logarithms are natural.

\subsection{The configuration model}

We recall the configuration model, which provides a convenient representation of uniform random regular graphs. In the $d$-regular configuration model, each vertex is assigned $d$ distinguishable half-edges, and a uniform perfect matching of the $dn$ half-edges is sampled. The resulting graph is a $d$-regular multigraph, possibly containing loops or multiple edges.

The probability that $t$ prescribed disjoint pairs of half-edges occur in the uniform matching is
\begin{equation}\label{eq:pairing-probability}
\prod_{j=0}^{t-1}\frac{1}{dn-2j-1}.
\end{equation}
For fixed $d$, the probability that the resulting multigraph is simple satisfies the following classical asymptotic estimate (see \cite{janson2009simple}):
\begin{equation}\label{eq:partitoShlein}
    \Pp(G^{\mathrm{CM}}_{n,d}\text{ is simple})
\longrightarrow p_d := \exp\left(-\frac{d^2-1}{4}\right).
\end{equation}

Moreover, conditional on simplicity, the resulting graph is uniformly distributed over $\Gnd$; see~\cite{Wormald1999}. Consequently, if an event $A_n$ holds with high probability in the configuration model, then
\[
\begin{aligned}
\Pp(G_{n,d}\notin A_n)
&=\Pp(G^{\mathrm{CM}}_{n,d}\notin A_n\mid G^{\mathrm{CM}}_{n,d}\text{ is simple})\\
&\leq
\frac{\Pp(G^{\mathrm{CM}}_{n,d}\notin A_n)}
{\Pp(G^{\mathrm{CM}}_{n,d}\text{ is simple})}
\longrightarrow0.
\end{aligned}
\]
Thus, high-probability statements transfer from the configuration model to the uniform simple model, although the converse need not hold.

\subsection{Contiguity of random graph models}

Contiguity allows us to transfer high-probability statements from auxiliary models with prescribed spanning structures to uniform random regular graphs. Two sequences of probability measures $(\mathcal P_n)$ and $(\mathcal Q_n)$ on the same spaces are
\emph{mutually contiguous}, written $\mathcal P_n\approx\mathcal Q_n$, if, for every sequence of events $A_n$,
\[
\mathcal P_n(A_n)\longrightarrow0
\qquad\Longleftrightarrow\qquad
\mathcal Q_n(A_n)\longrightarrow0.
\]
Thus, the two models share the same high-probability events. Mutual contiguity is also preserved under fixed finite products, allowing its application to independent pairs of graphs.

For the planted models considered below, if $Y_n(G)$ counts the admissible spanning structures in $G$, the unmarked planted distribution is the $Y_n$-size-biased version of the original law \cite{janson1995contiguity}, namely
\[
\frac{d\mathcal Q_n}{d\mathcal P_n}(G)
=\frac{Y_n(G)}{\Ee[Y_n]}.
\]
Controlling this normalised count via small-subgraph conditioning is a standard approach to establishing contiguity for random regular graphs; see~\cite{Wormald1999}.

This differs from conditioning the configuration model on simplicity: since simplicity has probability bounded away from zero, high-probability statements transfer to the uniform simple model, but mutual contiguity does not follow. In particular, this argument does not characterise the residual graph after selecting a spanning subgraph.

\subsection{Hamilton cycles and short cycle factors}

Our lower bounds rely on spanning structures that can be aligned between two independent graphs. The two constructions use these structures in different ways.  For Hamilton cycles, we shall work in a planted model and transfer the resulting graph event by contiguity. For logarithmic cycle factors, we only need an existence theorem; the residual pseudorandomness required later will be proved uniformly over all factors. The first input is the classical Hamiltonicity theorem of Robinson and Wormald~\cite{RobinsonWormald1994}.

\begin{theorem}[Hamiltonicity of random regular graphs]\label{thm:hamiltonicity}
For every fixed $d\geq3$, a uniform simple $d$-regular graph is Hamiltonian with high probability.
\end{theorem}

We shall also use the corresponding planted-model statement.  Let $H$ be a uniform Hamilton cycle on $[n]$ and let $R$ be an independent uniform simple $(d-2)$-regular graph.  Define the event
\[
D:=\{E(H)\cap E(R)=\varnothing\}.
\]
Conditionally on $D$, the union $H\cup R$ is a simple $d$-regular graph.  The following contiguity result is due to Wormald; see~\cite[Theorem~4.15(i)]{Wormald1999}.

\begin{theorem}[Hamilton-cycle contiguity]\label{thm:hamiltonian-contiguity}
Fix $d\geq3$ and let $n$ be admissible.  The distribution of
\(
H\cup R\mid D
\)
is mutually contiguous to $\mathcal U_{n,d}$.
\end{theorem}

The second spanning structure is a $k$-cycle factor, namely a spanning collection of $n/k$ vertex-disjoint cycles, each of length $k$. Kahn, Lubetzky, and Wormald \cite{KahnLubetzkyWormald2017} establish a sharp threshold for such factors in random cubic graphs, and state that the argument extends to every fixed degree $d\geq 3$. Before stating such an extension, we compute the degree-dependent first moment, which identifies the corresponding threshold location.

\begin{proposition}[First moment for $k$-cycle factors]\label{prop:cycle-factor-first-moment}
Fix $d\ge 3$, and let $k=k(n)\to\infty$ satisfy $k\mid n$ and $k=O(\log n)$. Put $q=n/k$, and let $X_{n,k}$ denote the number of \emph{realised} spanning $k$-cycle factors in the $d$-regular configuration model on $[n]$, where a realisation records the half-edge pairs used by the factor. Define
\[
 a_d:=\log(d-1)+\frac{d-2}{2}\log\!\left(\frac{d-2}{d}\right).
\]
Then
\begin{equation}\label{eq:cycle-factor-first-moment}
 \log \Ee X_{n,k}
 =n a_d-\frac nk\log(2n/e)+O_d(\log n).
\end{equation}
In particular, if $k\sim c\log n$ and $c<a_d^{-1}$, then $\Pp(X_{n,k}>0)\to0$.
\end{proposition}

\begin{proof}
There are exactly
\[
 |\mathcal F_{n,k}|=\frac{n!}{(2k)^q q!}
\]
labelled $k$-cycle factors on $[n]$. For a fixed factor $F$, choose at each vertex the two ordered half-edges that realise its two incident factor edges. This gives $[d(d-1)]^n$ half-edge realizations. Once the $n$ factor pairs have been selected, the remaining $(d-2)n$ half-edges may be paired arbitrarily. Therefore
\[
 \Ee X_{n,k}
 =\frac{n!}{(2k)^q q!}[d(d-1)]^n
 \frac{((d-2)n-1)!!}{(dn-1)!!}.
\]
Using $(2r-1)!!=(2r)!/(2^r r!)$ and Stirling's formula gives
\begin{align*}
 \log \Ee X_{n,k}
 &={}n\log[d(d-1)]+\log n!-q\log(2k)-\log q!\\
 &\quad+\log((d-2)n-1)!!-\log(dn-1)!!\\
 &={}n\left[\log(d-1)+\frac{d-2}{2}
       \log\!\left(\frac{d-2}{d}\right)\right]
 -\frac nk\log(2n/e)+O_d(\log n),
\end{align*}
which is \eqref{eq:cycle-factor-first-moment}.

If $k\sim c\log n$ with $c<a_d^{-1}$, then for some fixed $\eta>0$, $ka_d\le(1-\eta)\log n$ for all sufficiently large $n$. Hence $\log \Ee X_{n,k}\to-\infty$, and Markov's inequality gives $\Pp(X_{n,k}>0)\to0$ in the configuration model.  Since simplicity has probability bounded away from zero, the same non-existence conclusion holds after conditioning on simplicity.
\end{proof}

The preceding computation identifies the fixed-degree threshold constant. In particular, the natural location of the transition is
\begin{equation}\label{eq:Kd-def}
 K_d(n):=\frac{\log(2n/e)}{a_d}.
\end{equation}

\begin{theorem}[Short cycle factors in fixed degree]\label{thm:cycle-factor-existence}
Fix $d\ge3$ and let
\[
 a_d=\log(d-1)+\frac{d-2}{2}\log\!\left(\frac{d-2}{d}\right),
 \qquad
 K_d(n)=\frac{\log(2n/e)}{a_d}.
\]
Let $k=k(n)$ satisfy $k\mid n$,
\[
 k=K_d(n)+\omega(1),
 \qquad
 k=(1+o(1))K_d(n).
\]
Then a uniform simple $d$-regular graph contains a $k$-cycle factor with high probability. Conversely, for every fixed $\delta>0$, no $k$-cycle factor exists with high probability whenever
\[
 k\le \left(\frac1{a_d}-\delta\right)\log n.
\]
\end{theorem}

For $d=3$, the existence follows from \cite[Theorem~1]{KahnLubetzkyWormald2017}. Moreover, the authors state that their proof extends, with very few modifications, to every fixed degree $d\geq 3$, yielding the corresponding cycle-factor threshold and its two-point concentration. This is in line with the result we obtain in \cref{prop:cycle-factor-first-moment}, and we use the statement in \cite{KahnLubetzkyWormald2017} as a witness of the $k$-cycle factor existence with $k\geq 3$.

We note that the complementary non-existence result follows directly from \cref{prop:cycle-factor-first-moment} and Markov's inequality, independently of the methods used in~\cite{KahnLubetzkyWormald2017} to establish existence in cubic graphs.

\section{Low degrees and first-order overlap}\label{sec:easy_section}

The effect of optimising over vertex permutations depends strongly on the degree. We first settle the cases $d=1$ and $d=2$, where the component structure can be exploited directly. We then establish the first-order overlap for every fixed $d\geq3$. In this regime, Hamiltonicity gives a baseline of $n$ common edges, while a simple first-moment argument shows that the maximum overlap cannot be considerably larger at first order. The sharper upper and lower bounds on the surplus are developed in the following sections.

\subsection{Degrees one and two}

For $d=1$, every graph is a perfect matching. Since any two perfect matchings on $[n]$ are isomorphic, $M_n=n/2$ and $\GED_E(G_1,G_2)=0$ deterministically. For $d=2$, each graph is a disjoint union of cycles of length at least three. 

For a $2$-regular graph $G$, let $c(G)$ denote its number of cycles, and let $\gamma_\ell(G)$ denote its number of cycles of length $\ell$. The following estimate describes the cycle counts needed for our argument.

\begin{lemma}[Cycle counts in random $2$-regular graphs]\label{lem:degree-two-cycles}
Let $G_1,G_2$ be independent samples from $\mathcal U_{n,2}$. Then, for $i\in\{1,2\}$,
\[
c(G_i)=\left(\frac12+o_{\Pp}(1)\right)\log n.
\]
Moreover,
\[
\sum_{\ell=3}^n\min\{\gamma_\ell(G_1),\gamma_\ell(G_2)\}=O_{\Pp}(1).
\]
\end{lemma}

\begin{proof}
We first work in the $2$-regular configuration model. Write $\gamma$ for its total number of connected components, including loops and double edges. A uniform pairing can be constructed recursively by adding one vertex at a time: with probability $1/(2j-1)$, the new vertex forms a loop, and otherwise it is inserted into one of the edges of the pairing on the preceding $j-1$ vertices. Insertion preserves the number of components, whereas a loop creates one additional component. Consequently, $\gamma$ has the distribution of a sum of independent Bernoulli random variables with respective parameters $1/(2j-1)$, $1\leq j\leq n$. In particular,
\[
\E{\gamma}=\sum_{j=1}^n\frac1{2j-1}=\frac12\log n+O(1),\qquad \operatorname{Var}(\gamma)=O(\log n).
\]
Chebyshev's inequality implies $\gamma=(1/2+o_{\Pp}(1))\log n$. Since the configuration model is simple with probability bounded away from zero, the same conclusion holds for the number of cycles of a uniform simple $2$-regular graph.

To prove the second assertion, let $\gamma_\ell^{\mathrm{CM}}$ be the number of components of length $\ell$ in the configuration model. Counting cycles on $\ell$ prescribed vertices and using the pairing probability gives, for $1\leq\ell\leq n$,
\[
\E{\gamma_\ell^{\mathrm{CM}}}=\frac{1}{2\ell}\prod_{j=0}^{\ell-1}\frac{n-j}{n-j-\frac12}\leq\frac{K}{\ell}\sqrt{\frac{n}{n-\ell+1}},
\]
where $K$ is an absolute constant. Conditioning on simplicity increases the expectation of any non-negative random variable by at most the reciprocal of the simplicity probability. Hence, uniformly in $\ell$,
\[
\E{\gamma_\ell(G_i)}\leq\frac{K'}{\ell}\sqrt{\frac{n}{n-\ell+1}}
\]
for another absolute constant $K'$. By independence and $\min(x,y)\leq xy$ for non-negative integers $x,y$,
\[
\E{\sum_{\ell=3}^n\min\{\gamma_\ell(G_1),\gamma_\ell(G_2)\}}\leq\sum_{\ell=3}^n\E{\gamma_\ell(G_1)}\E{\gamma_\ell(G_2)}\leq (K')^2\sum_{\ell=3}^n\frac{n}{\ell^2(n-\ell+1)}=O(1).
\]
The last sum is bounded by splitting at $\ell=n/2$: the first part is $O(\sum_{\ell\geq1}\ell^{-2})$, and the second is $O(\log n/n)$. Markov's inequality completes the proof.
\end{proof}

We now prove \cref{thm:overlap-d2}. For the lower bound, we remove one edge from each cycle and concatenate the resulting paths to construct an alignment. For the upper bound, we associate a bipartite multigraph with the common graph, whose vertices represent the cycles of both graphs. We control the number of connected components through an elementary first-moment argument, which yields the leading-order asymptotics.

\begin{proof}[Proof of \cref{thm:overlap-d2}]

For the lower bound, remove one edge from each cycle of $G_i$, obtaining $c_i$ vertex-disjoint paths. Concatenate these paths to form an ordering of the vertices of $G_i$, and align the two orderings. Among the $n-1$ consecutive pairs in either ordering, exactly $c_i-1$ are not edges of $G_i$. Hence
\begin{equation}
M_n\ge n-c_1-c_2+1.
\label{eq:d2-lower}
\end{equation}

Next, we consider the upper bound. Fix a permutation $\pi$ and identify the vertices
through it. Let
\[
J=\pi(G_1)\cap G_2
\]
be the common graph. Construct a bipartite multigraph $B$ whose vertices represent the cycles of $\pi(G_1)$ and $G_2$, with one edge for each connected component of $J$,
including isolated vertices.

Let $b$ be the number of connected components of $B$, and let $r$ be the number of cyclic components of $J$. Since $J$ has maximum degree two, its number of
connected components is
\[
n-\operatorname{Ov}_\pi(G_1,G_2)+r.
\]
This is also the number of edges of $B$. As $B$ has $c_1+c_2$ vertices and $b$ connected
components, it follows that $n-\operatorname{Ov}_\pi(G_1,G_2)+r
\ge c_1+c_2-b$. Consequently,
\begin{equation}
\operatorname{Ov}_\pi(G_1,G_2)
\le n-c_1-c_2+b+r.
\label{eq:d2-deterministic}
\end{equation}
Since every cyclic component of $J$ must be an entire cycle of the same length in both graphs, from \cref{lem:degree-two-cycles} it follows $r = O_{\mathbb P}(1)$, uniformly in $\pi$. It remains to control $b$.

For $i=1,2$ and $0\le k\le n$, let $A_{i,k}$ be
the number of $k$-vertex subsets of $G_i$ that are
unions of cycles. Define
\[
T=\sum_{k=0}^n A_{1,k}A_{2,k}.
\]
Thus $T$ counts pairs of unions of cycles, one
from each graph, having the same number of vertices.

Every connected component of $B$ corresponds to
a nonempty set of vertices that is a union of cycles in both graphs. Taking arbitrary unions of
these $b$ disjoint sets gives $2^b$ distinct pairs
counted by $T$. In particular,
\begin{equation}
2^b\le T.
\label{eq:d2-components}
\end{equation}
We show that $\mathbb E[T]=O(\log n)$. First, consider the $2$-regular configuration model. Applying a standard argument for the configuration model, one obtains
\[
\mathbb E_{\mathrm{CM}}[A_{i,k}]
=\binom nk
\frac{(2k-1)!!(2(n-k)-1)!!}{(2n-1)!!}
\]
Using the standard estimate $\binom{2m}{m}\asymp 4^m/\sqrt{m+1}$, we obtain, uniformly in $0\le k\le n$,
\[
\mathbb E_{\mathrm{CM}}[A_{i,k}]
=\frac{\binom{2k}{k}\binom{2(n-k)}{n-k}}
{\binom{2n}{n}}
=O\left(\sqrt{\frac{n+1}{(k+1)(n-k+1)}}\right).
\]
Since the probability of simplicity is bounded away from zero, conditioning on simplicity increases these expectations by at most a constant factor $C$. Hence, using independence of $G_1,G_2$ we conclude that
\begin{align*}
\mathbb E[T]
&=\sum_{k=0}^n
  \mathbb E[A_{1,k}]\mathbb E[A_{2,k}]
\le C^2\sum_{k=0}^n
  \mathbb E_{\mathrm{CM}}[A_{1,k}]\mathbb E_{\mathrm{CM}}[A_{2,k}] = O(\log n).
\end{align*}
By Markov's inequality,
$T=O_{\mathbb P}(\log n)$, which, together with
\eqref{eq:d2-components}, yields
\[
b\le\log_2\log n+O_{\mathbb P}(1).
\]
This bound holds simultaneously
for all permutations $\pi$, since $T$ depends
only on the two original graphs. Combining this with \eqref{eq:d2-deterministic}
and $r = O_{\mathbb P}(1)$, we obtain
\[
M_n\le
n-c_1-c_2+O_{\mathbb P}(\log\log n).
\]
This proves the claim.\end{proof}

It is possible to obtain a stronger upper bound than the one derived in the proof above. In particular, $M_n \leq n-c_1-c_2+O_{\mathbb P}(1)$ . In fact, let
$\mathcal L_i$ be the set of subset sums of the cycle lengths of $G_i$, and put
\[
Z=
\big|(\mathcal L_1\cap\mathcal L_2)
\setminus\{0,n\}\big|.
\]

Each connected component of $B$ corresponds to a nonempty vertex set that is a union of entire cycles in both graphs. Ordering these $b$ sets and taking successive partial unions produces $b-1$ distinct nontrivial common subset sums. Hence
\begin{equation}\label{eq:bZ}
    b\le Z+1.
\end{equation}
The probabilistic estimate $Z=O_{\mathbb P}(1)$ is established in \cref{app:ewens-common-sums} of the appendix.

\subsection{Hamiltonian baseline and first-order overlap}

We now consider fixed $d\geq3$. By \cref{thm:hamiltonicity}, both $G_1$ and $G_2$ contain Hamilton cycles with high probability. Mapping one Hamilton cycle onto the other preserves all $n$ cycle edges and therefore gives $M_n\geq n$ with high probability. This is the baseline that we shall improve in Sections~5 and~6.

To establish the first-order result, it remains to show that $M_n$ cannot exceed $(1+\eps)n$ with high probability. A simple union bound suffices; Section~4 will refine this argument to capture the surplus of order $n/\log n$.

\begin{proof}[Proof of \cref{thm:first-order-overlap}]
The lower bound follows from the Hamiltonian alignment described above. For the upper bound, fix $\eps>0$ sufficiently small that $2(1+\eps)<d$, and set $t=\lceil(1+\eps)n\rceil$. Condition on $G_1$. If $M_n\geq t$, there exist a permutation $\pi\in\Sn$ and a set $F\subseteq E(G_1)$ of $t$ edges such that $\pi(F)\subseteq E(G_2)$.

For a prescribed set of $t$ edges, the configuration model gives at most $d^{2t}$ choices of half-edges to realise their incidences. By \eqref{eq:pairing-probability}, the probability of any prescribed pairing is at most
\[
\prod_{j=0}^{t-1}\frac{1}{dn-2j-1}\leq\left(\frac{1}{dn-2t+1}\right)^t.
\]
Since $dn-2t$ is bounded below by a positive constant times $n$, and conditioning on simplicity costs only a bounded factor, there is a constant $K=K(d,\eps)$ such that
\[
\Pp(\pi(F)\subseteq E(G_2))\leq\left(\frac Kn\right)^t.
\]
A union bound over the $n!$ permutations and the $\binom{nd/2}{t}$ choices of $F$ now yields
\[
\Pp(M_n\geq t)\leq n!\binom{nd/2}{t}\left(\frac Kn\right)^t.
\]
Using $\binom{nd/2}{t}\leq2^{nd/2}$ and Stirling's formula, we obtain
\[
\log\Pp(M_n\geq t)\leq n\log n-t\log n+O_{d,\eps}(n)=-\eps n\log n+O_{d,\eps}(n)\longrightarrow-\infty.
\]
Thus $\Pp(M_n\geq(1+\eps)n)\to0$. Since it suffices to prove this for arbitrarily small $\eps>0$, the Hamiltonian lower bound implies $M_n/n\xrightarrow{\Pp}1$.
\end{proof}

\section{The upper bound}\label{sec:upper}

We now refine the first-moment argument of \cref{sec:easy_section} to the scale of the surplus above the Hamiltonian baseline.  The first-order argument counts arbitrary common edge sets.  At second order, this is too wasteful: if an aligned intersection contains $n+s$ edges, then the relevant quantity is not its total size but its graph excess.  We therefore first remove the tree-like part of the intersection and reduce the problem to counting common subgraphs of minimum degree at least two. For a finite graph $H$, write
\[
    \delta(H):= \min_{u \in V(H)} \text{deg}(u), \qquad \Delta(H):= \max_{u \in V(H)} \text{deg}(u), \qquad \operatorname{ex}(H):=|E(H)|-|V(H)|.
\]
The proof of \cref{thm:second-order-upper} has three steps.  Suppose that $M_n\ge n+s$.
\begin{enumerate}
    \item \emph{Extract a structural witness.}  From an aligned intersection with at least $n+s$ edges, we extract a common graph $H$ (a core) with minimum degree at least two and $\operatorname{ex}(H)=s$. Thus, the surplus becomes exactly the excess of a common core.

    \item \emph{Count common cores.}  If $H$ has $v$ vertices, let $Z_{v,s}$ count such common core witnesses. A configuration-model enumeration gives, uniformly for $s=O(n/\log n)$,
    \[
        \log \Ee Z_{v,s}
        \le n f_d(v/n)-s\log n
        +O_d\!\left(s\log\frac ns+\log n\right),
    \]
    where the leading term $n f_d(v/n)$ quantifies the entropy of a core occupying a fraction $v/n$ of the vertices.

    \item \emph{Optimise over the core size.} The entropy function $f_d$ is maximised at $v/n=d/(d+1)$, attaining value $b_d$. Markov's inequality and a union bound over $v$ conclude the proof.
\end{enumerate}

\subsection{From surplus overlap to a common core}
We begin by extracting a common core between $G_1$ and $G_2$ with prescribed excess.

\begin{lemma}[Exact-excess core reduction]\label{lem:core-reduction}
Let $J$ be a simple graph on at most $n$ vertices with at least $n+s$ edges, where $s\ge1$, and suppose $\Delta(J)\le d$.  Then $J$ contains a nonempty subgraph $H$ such that
\[
    \delta(H)\ge2,
    \qquad
    \Delta(H)\le d,
    \qquad
    \operatorname{ex}(H)=s.
\]
\end{lemma}

\begin{proof}
Add isolated vertices, if necessary, so that $J$ is viewed as a graph on exactly $n$ vertices.  Repeatedly delete vertices of degree zero or one, together with all incident edges.  Deleting an isolated vertex increases $|E|-|V|$ by one, while deleting a degree-one vertex leaves $|E|-|V|$ unchanged.  Hence the resulting $2$-core $K$ satisfies
\[
    \delta(K)\ge2,
    \qquad
    \operatorname{ex}(K)\ge \operatorname{ex}(J)\ge s.
\]
Since $s>0$, the core is nonempty.

If $\operatorname{ex}(K)=s$, take $H=K$.  Otherwise, choose a connected component $C$ of $K$ with positive excess.  Such a component contains a non-bridge edge $e$. Delete $e$ and then repeatedly prune every degree-one vertex that is created.  Deleting $e$ lowers the excess by exactly one, while each subsequent pruning step removes one vertex and one edge and therefore leaves the excess unchanged.  Moreover, $C-e$ still contains a cycle, so the pruning procedure leaves a nonempty graph of minimum degree at least two. Thus, one application of this operation lowers the total excess by exactly one and never increases the maximum degree. Iterating until the total excess is $s$ returns $H$.
\end{proof}

\begin{corollary}\label{cor:overlap-core}
If $M_n\ge n+s$, then for some $1\le v\le n$ there exist vertex sets $V_1,V_2\subseteq[n]$ of size $v$, a bijection $\phi:V_1\to V_2$, and a simple graph $H$ on $V_1$ such that
\[
    |E(H)|=v+s,
    \qquad
    2\le \text{deg}_H(u)\le d\quad (u\in V_1),
\]
and
\[
    H\subseteq G_1[V_1],
    \qquad
    \phi(H)\subseteq G_2[V_2].
\]
\end{corollary}

\begin{proof}
Choose a permutation attaining at least $n+s$ common edges and identify the corresponding aligned intersection with a graph $J$ on $[n]$.  Since both ambient graphs are $d$-regular, $\Delta(J)\le d$.  Apply \cref{lem:core-reduction} and then undo the identification of the two vertex sets.
\end{proof}

The point of \cref{cor:overlap-core} is that the surplus is now encoded by a much more rigid witness.  Tree-like extensions of the aligned intersection disappear under the $2$-core reduction, while the quantity $|E|-|V|$, which records the surplus, is preserved.

\subsection{A first-moment bound for common cores}

We first work with two independent $d$-regular configuration models. Conditioning both models on simplicity changes the upper bounds on the probability only by a constant factor $p_d$ defined in \eqref{eq:partitoShlein}.

Fix integers $v\le n$ and $s\ge1$, and put $m=v+s$.  Let $Z_{v,s}$ count choices of vertex sets $V_1,V_2\subseteq[n]$ of size $v$, a bijection $\phi:V_1\to V_2$, and a simple graph $H$ on $V_1$ having $m$ edges and all degrees in $\{2,\ldots,d\}$, such that
\[
    H\subseteq G_1[V_1],
    \qquad
    \phi(H)\subseteq G_2[V_2].
\]

Define, for $0\le x\le1$,
\begin{equation}\label{eq:f-def}
    f_d(x)
    :=-2(1-x)\log(1-x)
      +2x\log[d(d-1)]
      -d\log d
      +(d-2x)\log(d-2x),
\end{equation}
with the convention $0\log0=0$. The following proposition contains all the technical enumeration needed for the upper bound.

\begin{proposition}[First moment for common cores]\label{prop:core-first-moment}
Fix $d\ge3$.  Uniformly over feasible integers $v$ and over integers $s$ satisfying $1\le s=O(n/\log n)$,
\begin{equation}\label{eq:core-first-moment}
    \log \Ee Z_{v,s}
    \le
    n f_d(v/n)-s\log n
    +O_d\!\left(s\log\frac ns+\log n\right).
\end{equation}
\end{proposition}

The proof is not particularly complex, but rather technical, and it is given in \cref{app:upper-bound-count}. We briefly indicate where the three terms in \eqref{eq:core-first-moment} come from.  The quantity $n f_d(v/n)$ is the exponential entropy of an excess-zero core occupying a fraction $v/n$ of the vertices. Raising the excess from $0$ to $s$ changes the pairing probability by a factor $n^{-s}e^{O_d(s)}$, producing the term $-s\log n$.  Finally, the $2s$ units of degree above the minimum value two can be distributed among the $v$ vertices in at most
\[
    \exp\!\left\{O_d\!\left(s\log\frac ns\right)\right\}
\]
ways. Thus, the main exponential competition is between the entropy of the core size and the $n^{-s}$ cost of the surplus.

\subsection{Optimising the core size}

The optimisation left by \cref{prop:core-first-moment} is elementary but important: it identifies which core sizes are most ``dangerous'' for the first-moment argument.

\begin{lemma}[Optimising the core entropy]\label{lem:core-opt}
The function $f_d$ is strictly concave on $[0,1]$ and has unique maximiser
\begin{equation}\label{eq:bd-def}
    x_d^*=\frac d{d+1},\qquad f_d(x_d^*)
    =d\log(d-1)-(d-2)\log(d+1)
    =:b_d.
\end{equation}
\end{lemma}

\begin{proof}
Differentiating \eqref{eq:f-def},
\[
    f_d'(x)
    =2\log\!\left(
      \frac{d(d-1)(1-x)}{d-2x}
    \right),
\]
and
\[
    f_d''(x)
    =-\frac{2(d-2)}{(1-x)(d-2x)}<0,
    \qquad 0<x<1.
\]
Thus $f_d$ is strictly concave.  Solving $f_d'(x)=0$ gives $x=d/(d+1)$, and substituting this value into \eqref{eq:f-def} yields \eqref{eq:bd-def}.
\end{proof}

\begin{theorem}[Second-order upper tail]\label{thm:second-order}
Fix $d\ge3$ and let $1\le s=O(n/\log n)$ be an integer.  Then
\begin{equation}\label{eq:second-order-tail}
    \log\Pp(M_n\ge n+s)
    \le
    b_dn-s\log n
    +O_d\!\left(s\log\frac ns+\log n\right).
\end{equation}
\end{theorem}

\begin{proof}
By \cref{cor:overlap-core}, on the event $\{M_n\ge n+s\}$ there is at least one common core counted by $Z_{v,s}$ for some $v$. \cref{prop:core-first-moment} and \cref{lem:core-opt} imply, uniformly over feasible $v$,
\[
    \log \Ee Z_{v,s}
    \le
    b_dn-s\log n
    +O_d\!\left(s\log\frac ns+\log n\right).
\]
There are at most $n$ choices for $v$.  Hence Markov's inequality and a union bound give \eqref{eq:second-order-tail}, up to an additional $O(\log n)$ term.  Conditioning both configurations on simplicity multiplies the probability by at most $(p_d+o(1))^{-2}$, which is absorbed in the error term.
\end{proof}

\begin{proof}[Proof of \cref{thm:second-order-upper}]
Apply \cref{thm:second-order}. Next, take
\[
    s=\left\lceil(b_d+\varepsilon)\frac n{\log n}\right\rceil.
\]
Since
\[
    s\log\frac ns
    =O\!\left(\frac{n\log\log n}{\log n}\right)
    =o(n),
\]
we obtain
\[
    \log\Pp(M_n\ge n+s)
    \le -\varepsilon n+o(n)\longrightarrow-\infty.
\]
\end{proof}
Recall that the lower-bound constant appearing in \cref{thm:second-order-lower} is
\[
    a_d
    :=\log(d-1)+\frac{d-2}{2}\log\!\left(\frac{d-2}{d}\right).
\]
Evaluating the entropy function at the spanning value $x=1$ gives
\[
    f_d(1)=2a_d.
\]
If one knew that every overlap with surplus $s$ contained a spanning witness with minimum degree two, the same first-moment calculation would lead to the smaller upper coefficient $2a_d$.

The core reduction gives no such information: pruning degree-zero and degree-one vertices may remove a linear fraction of the vertices.  Once all core sizes are allowed, the entropy function is maximised instead at
\[
    x_d^*=\frac d{d+1}<1,
\]
and the resulting coefficient is
\begin{equation}
    b_d=f_d(x_d^*)>f_d(1)=2a_d.\label{remark:2ad}
\end{equation}

Indeed,
\[
    b_d-2a_d
    =(d-2)\log\!\left(
      \frac{d(d-1)}{(d+1)(d-2)}
    \right)>0.
\]
The interpretation is only a statement about the present first-moment bound: among the witnesses over which we union-bound, non-spanning cores of size about $dn/(d+1)$ carry the largest entropy.  It does not prove that the $2$-core of an optimal alignment typically has this size. Improving $b_d$ toward $2a_d$, or below it, would therefore require additional structural information that rules out or more sharply penalises these non-spanning witnesses.
\section{Alignment of Hamiltonian cycles}
\label{sec:hamiltonian-alignment}

In \cref{sec:easy_section}, Hamiltonicity gave the baseline $M_n\ge n$ with high probability.  We now ask how much one can improve this overlap while keeping the Hamiltonian alignment essentially fixed.  More precisely, once Hamilton cycles
\[
H_1=(v_0,\ldots,v_{n-1})\subseteq G_1,
\qquad
H_2=(w_0,\ldots,w_{n-1})\subseteq G_2
\]
have been chosen, we only allow cyclic shifts of the identification between them.  For $s\in\mathbb Z_n$, write
\[
\rho_s(x)=x+s\pmod n,
\]
and use the same notation for the induced action on edges and edge sets.  In cyclic coordinates, the corresponding alignment is
\[
\pi_s(v_i)=w_{i+s},\qquad i\in\mathbb Z_n.
\]
Every such shift preserves all $n$ Hamilton-cycle edges. The first result shows that these $n$ restricted alignments already yield a substantial surplus.

\begin{theorem}[Hamilton-cycle surplus]
\label{thm:hamiltonian-rotation-overlap}
Fix $d\ge3$ and $\varepsilon\in(0,1)$. Let $G_1,G_2$ be independent uniform simple $d$-regular graphs on $[n]$, with $n$ even. With high probability, there exist Hamilton cycles
\[
H_1=(v_0,\ldots,v_{n-1})\subseteq G_1,
\qquad
H_2=(w_0,\ldots,w_{n-1})\subseteq G_2,
\]
and a shift $s\in\mathbb Z_n$ such that
\[
|\pi_s(E(G_1))\cap E(G_2)|
\ge
n+\left\lfloor(1-\varepsilon)\frac{\log n}{\log\log n}\right\rfloor.
\]
\end{theorem}
In fact, the planted-model argument gives more than existence: a typical pair of Hamilton cycles works.

\begin{theorem}[Typical Hamilton-cycle pairs]
\label{thm:typical-hamiltonian-alignment}
Fix $d\ge3$ and $\varepsilon\in(0,1)$, and let $G_1,G_2$ be independent uniform simple $d$-regular graphs on $[n]$, with $n$ even. With high probability, all but an $o(1)$ proportion of pairs of Hamilton cycles
\[
H_1\subseteq G_1,
\qquad
H_2\subseteq G_2,
\]
admit a cyclic alignment whose overlap is at least
\[
n+\left\lfloor(1-\varepsilon)\frac{\log n}{\log\log n}\right\rfloor.
\]
Equivalently, after sampling $ G_1$ and $ G_2$, if one Hamilton cycle is then chosen uniformly from each graph, independently conditional on the graphs, the conditional probability that the selected pair admits such an alignment tends to one in probability.
\end{theorem}
The precise formulation in terms of the random proportion of good Hamilton-cycle pairs, together with the proof, is given in \cref{app:hamiltonian-technical}. The argument to prove \cref{thm:hamiltonian-rotation-overlap} makes use of three conceptual steps.
\begin{enumerate}
    \item \emph{Reduction to a planted model.} We first work with a contiguous random graph defined through
    \[
    G=H\cup R,
    \]
    where $H$ is a uniform Hamilton cycle and $R$ is a uniform simple $(d-2)$-regular graph.

    \item \emph{Rotations of the Hamilton cycle.} Once the two planted Hamilton cycles are identified with the standard cycle on $\mathbb Z_n$, every cyclic shift preserves their $n$ edges. The problem, therefore, reduces to finding a shift $s$ for which
    \[
    |\rho_s(E(R_1))\cap E(R_2)|
    \]
    is large. From each residual graph, we extract a uniformly distributed matching, and the main probabilistic input is provided by showing that among their $n$ cyclic shifts, one has overlap of order $\log n/\log\log n$.

    \item \emph{Returning to $\mathcal U_{n,d}$.} The residual-overlap event survives the edge-disjointness conditioning because that conditioning has probability bounded away from zero.  We then forget the distinguished Hamilton cycles and transfer the resulting unmarked graph event to the uniform model by contiguity.
\end{enumerate}

\subsection{The random-matching input}
\label{subsec:ham-random-matching}

\begin{lemma}[A successful cyclic shift]
\label{lem:cyclic-matching-witness}
Fix $c\in(0,1/2]$ and $\varepsilon\in(0,1)$. For $i=1,2$, let $J_i$ be independent uniform matchings of sizes $a_i$ on $\mathbb Z_n$, where
\[
cn\le a_i\le n/2.
\]
Set
\[
k=\left\lfloor(1-\varepsilon)\frac{\log n}{\log\log n}\right\rfloor.
\]
Then
\[
\Prob{\max_{s\in\mathbb Z_n}|\rho_s(J_1)\cap J_2|\ge k}\longrightarrow1,
\]
uniformly over the allowed values of $a_1,a_2$.
\end{lemma}

Let us briefly explain the mechanism behind the lemma. For an edge $e=\{u,v\}$, define its cyclic length by
\[
\ell(e)=\min\{(v-u)\bmod n,(u-v)\bmod n\}.
\]
Rotations preserve cyclic length. From $J_1$, one can retain a linear-sized submatching $F$ containing at most one edge of each non-antipodal cyclic length.  Consequently, the translations
\[
F_s=\rho_s(F),\qquad s\in\mathbb Z_n,
\]
are pairwise edge-disjoint.  For a fixed shift $s$, a $k$-subset of $F_s$ that is contained in $J_2$ is a certificate that the shift has overlap at least $k$.  Counting these certificates over all shifts and applying a second-moment argument proves \cref{lem:cyclic-matching-witness}. The distinct-length reduction is precisely what makes the dependence between different shifts manageable; see \cref{app:hamiltonian-technical} for the details.

\subsection{The planted Hamilton-cycle model}
\label{subsec:ham-planted}

Set $q=d-2$. For each $i=1,2$, independently sample a uniform Hamilton cycle $H_i$ on $[n]$ and an independent uniform simple $q$-regular graph $R_i$ on $[n]$. Define
\[
D_i:=\{E(H_i)\cap E(R_i)=\varnothing\}.
\]
Conditionally on $D_i$, the union
\[
G_i=H_i\cup R_i
\]
is a simple $d$-regular graph.  Denote by $\widehat{\mathcal Q}_{n,d}$ the law of the marked pair $(G_i,H_i)$ after conditioning on $D_i$, and by $\mathcal Q_{n,d}$ its unmarked graph marginal. By \cref{thm:hamiltonian-contiguity},
\[
\mathcal Q_{n,d}\approx\mathcal U_{n,d}.
\]

We now prove the main existence result.

\begin{proof}[Proof of \cref{thm:hamiltonian-rotation-overlap}]
Let
\[
k=\left\lfloor(1-\varepsilon)\frac{\log n}{\log\log n}\right\rfloor.
\]
We first work before conditioning on $D_1,D_2$.  Choose a cyclic ordering of each $H_i$, using only $H_i$, and identify its vertices with $\mathbb Z_n$.  Conditional on $H_1,H_2$, the residual graphs $R_1,R_2$ remain independent uniform simple $q$-regular graphs, since they are independent of the planted cycles and their laws are invariant under vertex permutations.

By the uniform matching-extraction lemma in \cref{app:hamiltonian-technical}, from each $R_i$ we may choose a matching
\[
J_i\subseteq E(R_i),\qquad |J_i|=\lfloor n/4\rfloor,
\]
such that $J_i$ is marginally uniform among all matchings of that size.  The two extracted matchings can be chosen independently. \cref{lem:cyclic-matching-witness} therefore yields
\[
\Prob{\max_{s\in\mathbb Z_n}|\rho_s(J_1)\cap J_2|\ge k}\longrightarrow1.
\]
Since $J_i\subseteq E(R_i)$,
\begin{equation}
\label{eq:residual-shift-success}
\Prob{\max_{s\in\mathbb Z_n}|\rho_s(E(R_1))\cap E(R_2)|\ge k}\longrightarrow1.
\end{equation}

The edge-disjointness lemma in \cref{app:hamiltonian-technical} gives, uniformly in the planted Hamilton cycles,
\[
\Prob{D_i}\longrightarrow e^{-q},
\]
and hence
\[
\Prob{D_1\cap D_2}\longrightarrow e^{-2q}>0.
\]
If $B$ denotes the failure event in \eqref{eq:residual-shift-success}, then
\[
\Prob{B\mid D_1\cap D_2}
\le
\frac{\Prob{B}}{\Prob{D_1\cap D_2}}
=o(1).
\]
Thus the residual surplus survives the conditioning.

Under $D_1\cap D_2$, identify both distinguished Hamilton cycles with
\[
C_n=\{\{x,x+1\}:x\in\mathbb Z_n\}.
\]
Every rotation preserves $C_n$, and the residual edges are disjoint from $C_n$.  Therefore, for every $s$,
\[
|\rho_s(E(G_1))\cap E(G_2)|
=n+|\rho_s(E(R_1))\cap E(R_2)|.
\]
Consequently, under the marked planted law $\widehat{\mathcal Q}_{n,d}^{\otimes2}$, with probability $1-o(1)$ the distinguished pair of Hamilton cycles admits a cyclic alignment with overlap at least $n+k$.

It remains to forget the distinguished cycles.  Consider the unmarked failure event that no pair of Hamilton cycles in $G_1, G_2$ admits a cyclic alignment with overlap at least $n+k$.  Under $\mathcal Q_{n,d}^{\otimes2}$ this event has probability $o(1)$, since the distinguished planted pair already succeeds with probability $1-o(1)$.  Mutual contiguity, preserved under finite products, now gives the same conclusion under $\mathcal U_{n,d}^{\otimes2}$.  This proves the theorem.
\end{proof}

The stronger \cref{thm:typical-hamiltonian-alignment} follows from the same marked planted law.  The additional observation is that, conditional on the unmarked graph $G$, the distinguished planted Hamilton cycle is uniform among all Hamilton cycles of $G$.  This turns the failure probability for the distinguished pair into the conditional proportion of bad Hamilton-cycle pairs.  We give the short formal argument in \cref{app:hamiltonian-technical}.

\subsection{Sharpness of the cyclic-shift scale}
\label{subsec:ham-sharpness}

For completeness, we record that the scale in \cref{lem:cyclic-matching-witness} is sharp at first order.  This is not needed for \cref{thm:hamiltonian-rotation-overlap}, but it identifies the natural limitation of the random-matching mechanism used above.

\begin{proposition}[Cyclic overlap of random matchings]
\label{prop:cyclic-matching-overlap}
Fix $c\in(0,1/2]$. Let $J_1,J_2$ be independent uniform matchings of sizes $a_1,a_2$ satisfying
\[
cn\le a_i\le n/2.
\]
With
\[
T_n:=\max_{s\in\mathbb Z_n}|\rho_s(J_1)\cap J_2|,
\]
we have
\[
\frac{T_n}{\log n/\log\log n}\xrightarrow{\mathbb P}1,
\]
uniformly over the allowed matching sizes.
\end{proposition}

The lower bound is \cref{lem:cyclic-matching-witness}; the matching upper tail is a direct first-moment estimate.  We defer the proof to \cref{app:hamiltonian-technical}.  The same estimate applies to unions of a fixed number of independent perfect matchings and yields \cref{thm:cyclic-matching-unions}; its proof is recorded there as well.
\section{A lower bound from cycle factors}\label{sec:lower}

The Hamilton-cycle argument of \cref{sec:hamiltonian-alignment} starts from a spanning cycle and uses contiguity to expose a tractable random residual graph.  In principle, one could proceed similarly with $k$-cycle factors: Kahn--Lubetzky--Wormald prove contiguity for the corresponding planted model in the cubic case \cite{KahnLubetzkyWormald2017}, and state that their argument extends to every fixed degree $d\ge3$.  Such an approach would allow one to study the residual graph of a particular factor and then transfer the resulting graph event to the uniform model.

Here, we use a different mechanism. Near the cycle factor existence threshold, the family of relevant factors is small enough that the residual properties needed for the alignment can be proved simultaneously for every factor contained in the graph. More precisely, the expected number of realised factors is only $\exp(o(n))$, whereas for a fixed realised factor the relevant residual-density estimate fails with probability $\exp(-\Theta(n))$.  This makes it possible to apply a union-bound over the entire family and then choose a factor after observing the graph. Recall
\[
a_d:=\log(d-1)+\frac{d-2}{2}\log\!\left(\frac{d-2}{d}\right).
\]
We first prove the result under $k\mid n$, and later remove the divisibility condition.

\begin{theorem}[Cycle-factor surplus]\label{thm:cycle-factor-overlap}
Fix $d\ge3$ and $\varepsilon>0$.  Let $k=k(n)$ satisfy
\begin{equation}\label{eq:k-regime}
k\mid n,
\qquad
k=\frac{\log n}{a_d}+o(\log n),
\qquad
k\ge \frac{\log n}{a_d}+\sqrt{\log n}.
\end{equation}
If $G_1,G_2$ are independent uniform simple $d$-regular graphs on $[n]$, then
\[
\Pp\left(
M_n \ge n+(1-\varepsilon)\frac nk
\right)\longrightarrow1.
\]
\end{theorem}
\noindent The proof of the latter theorem follows three conceptual steps.
\begin{enumerate}
    \item \emph{Control factors simultaneously.} \cref{prop:cycle-factor-first-moment} gives only $\exp(o(n))$ realized factors in expectation at the scale~\eqref{eq:k-regime}. Conditional on one realised factor, the unused half-edges form a uniform $(d-2)$-regular pairing.  An exponentially small concentration bound can therefore be union-bounded over all factors and all large unions of their cycles.

    \item \emph{Build common stars.}  The resulting residual density guarantees a factor cycle with linearly many attachment positions. A suitable rotation aligns many attachment positions of two centre cycles, while a local-sparsity estimate lets us choose $t$ distinct leaves on both sides.  A star using $t+1$ factor cycles preserves all factor edges and gains $t$ additional common residual edges.

    \item \emph{Iterate.}  Repeating the star construction while a positive fraction of factor cycles remains gives asymptotically one additional common edge per factor cycle.  Since $q=(a_d+o(1))n/\log n$, this yields the required surplus in the divisible case.
\end{enumerate}

\subsection{Uniform control over all short cycle factors}\label{subsec:uniform-cycle-factor-control}

Put $q:=\frac nk$. In the configuration model, a realised $k$-cycle factor records both the abstract labelled factor and the particular half-edge pairs realising its $n$ edges.  If $Y$ denotes the number of realized factors, \cref{prop:cycle-factor-first-moment} gives
\begin{equation}\label{eq:subexp-number-factors}
\log \mathbb E[Y]
=na_d-\frac nk\log(2n/e)+O_d(\log n)
=o(n),
\end{equation}
where the last equality follows from~\eqref{eq:k-regime}.  Thus
\[
\mathbb E[Y]=\exp(o(n)).
\]
This subexponential count is the key feature that allows us to work uniformly over all realised factors, rather than with a distinguished planted factor. Let $F$ be a $k$-cycle factor of a simple $d$-regular graph $G$, and write
\[
R:=G\setminus E(F)
\]
for the graph containing its residual edges.  We shall use the following two properties.

\medskip
\noindent\textbf{Property A (local sparsity).}
There is a constant $J=J_d$ such that every vertex set $S$ with $|S|\le2k$ spans at most $|S|+J$ edges in $G$.

\medskip
\noindent\textbf{Property B$(\eta)$ (residual density on large unions).}
For every union $U$ of cycles of $F$ with $|U|\ge\eta n$,
\[
e_R(U)\ge\frac{(d-2)|U|^2}{4n}.
\]

Property A has two consequences that will be used repeatedly.  A single factor cycle contains at most $J$ residual edges internally, and any two distinct factor cycles have at most $J$ residual edges between them: in the latter case, apply Property A to their union, which has $2k$ vertices and already contains exactly $2k$ factor edges.

\begin{proposition}[Uniformly good short cycle factors]\label{prop:all-cycle-factors-good}
Fix $d\ge3$ and $\eta>0$.  Under~\eqref{eq:k-regime}, with probability $1-o(1)$ a uniform simple $d$-regular graph satisfies Property A and every $k$-cycle factor it contains satisfies Property B$(\eta)$.
\end{proposition}

\begin{proof}
The local-sparsity statement is proved in \cref{lem:cf-local-sparsity} of \cref{app:cycle-factor-technical}.

For Property B, fix a particular realised factor $F$ in the configuration model.  Conditional on the occurrence of its $n$ prescribed factor pairs, the remaining $(d-2)n$ half-edges are uniformly paired.  Hence, for a fixed union $U$ of factor cycles with $|U|\ge\eta n$, \cref{lem:cf-pairing-concentration} gives
\[
\Pp\bigl(U\text{ violates Property B}(\eta)\mid F\text{ occurs}\bigr)
\le e^{-c_{d,\eta}n}
\]
for some constant $c_{d,\eta}>0$. The factor has $q$ cycles and therefore at most $2^q$ unions of cycles. Consequently,
\[
\Pp\bigl(F\text{ is bad}\mid F\text{ occurs}\bigr)
\le 2^q e^{-c_{d,\eta}n}.
\]
If $Y_{\rm bad}$ denotes the number of bad realised factors, summing this conditional estimate over all possible realised factors gives
\[
\mathbb E[Y_{\rm bad}]
\le
(\mathbb E[ Y])2^q e^{-c_{d,\eta}n}.
\]
By~\eqref{eq:subexp-number-factors}, $\mathbb E[Y]=e^{o(n)}$, while $2^q=e^{o(n)}$ because $q=O(n/\log n)$.  Hence
\[
\mathbb E[Y_{\rm bad}]\le e^{-c_{d,\eta}n+o(n)}=o(1).
\]
Markov's inequality shows that, with probability $1-o(1)$ in the configuration model, no realised factor is bad.  Conditioning on simplicity preserves this high-probability statement, since the probability of simplicity tends to a positive constant.  Combining it with Property A proves the proposition.
\end{proof}

The point of \cref{prop:all-cycle-factors-good} is that factor existence and factor regularity are now separated.  \cref{thm:cycle-factor-existence} ensures that a factor exists, while the proposition says that whichever factor is selected afterwards already has the residual properties required by the alignment argument.

\subsection{The deterministic common-star construction}\label{subsec:common-star}

We now forget randomness.  Let $G_1,G_2$ be two simple $d$-regular graphs on $n$ vertices, and suppose that each contains a $k$-cycle factor $F_i$ satisfying Properties A and B$(\eta)$ with the same constants $J$ and $\eta$.  Call a factor cycle \emph{available} if none of its vertices has yet been used.  At every stage, the two graphs have the same number of available cycles.

\begin{lemma}[A center with many attachment vertices]\label{lem:center-cycle}
Suppose that $m\ge\eta q$ cycles are available in one of the graphs.  For all sufficiently large $k$, some available cycle contains at least $\eta k/4$ distinct vertices having a residual edge to another available cycle.
\end{lemma}

\begin{proof}
Let $U$ be the union of the $m$ available cycles.  Then $|U|=mk\ge\eta n$.  Property B gives
\[
e_{G\setminus F}(U)
\ge\frac{(d-2)|U|^2}{4n}.
\]
By Property A, at most $Jm$ residual edges lie internally inside the individual factor cycles. Hence, the number $E_{\rm cross}$ of residual edges joining two distinct available cycles satisfies, for all sufficiently large $k$,
\[
E_{\rm cross}
\ge
\frac{(d-2)|U|^2}{4n}-Jm
\ge
\frac{(d-2)|U|^2}{8n}.
\]
Summing the crossing degrees over the $m$ available cycles counts every crossing edge twice. Therefore, some cycle is incident to at least
\[
\frac{2E_{\rm cross}}{m}
\ge
\frac{(d-2)|U|k}{4n}
\ge
\frac{(d-2)\eta k}{4}
\]
crossing residual edges. Each vertex has a residual degree at most $d-2$, so at least $\eta k/4$ distinct vertices of this cycle have a residual edge to another available cycle.
\end{proof}

\begin{lemma}[Common star]\label{lem:common-star}
Fix $\eta>0$ and a positive integer $t$.  For all sufficiently large $k$, whenever both graphs have at least $\eta q$ available cycles, one can select $t+1$ available cycles in each graph and a bijection between their vertex sets that preserves all $(t+1)k$ factor edges and at least $t$ additional residual edges, one from a designated centre cycle to each of $t$ distinct leaf cycles.
\end{lemma}

\begin{proof}
By \cref{lem:center-cycle}, choose a centre cycle $C_i$ in each graph such that at least $\eta k/4$ of its vertices have a residual edge to another available cycle.  Number the vertices of each centre cyclically by $\mathbb Z/k\mathbb Z$, and let $A_i$ be the corresponding sets of attachment positions.  Then
\[
|A_1|,|A_2|\ge\frac{\eta k}{4}.
\]
Averaging over the $k$ rotations gives
\[
\sum_{s\in\mathbb Z/k\mathbb Z}|A_1\cap(A_2+s)|
=|A_1||A_2|.
\]
Hence, some rotation of the identification of the two centres produces a set $B$ of simultaneously good positions satisfying
\begin{equation}\label{eq:aligned-positions}
|B|\ge\frac{|A_1||A_2|}{k}
\ge\frac{\eta^2k}{16}.
\end{equation}

We now inductively choose $t$ positions and $t$ distinct leaf cycles. Suppose $r<t$ leaves have already been selected in each graph. In graph $i$, call a position of $B$ forbidden if its centre vertex is joined by a residual edge to one of those previously selected leaves.  Property A bounds the number of residual edges between the centre and any fixed leaf by $J$. Thus, at most $rJ$ positions are forbidden in either graph, and at most $2rJ$ positions are forbidden in at least one of them.  Since $t$ is fixed, for sufficiently large $k$,
\[
\frac{\eta^2k}{16}>2(t-1)J.
\]
Therefore, a non-forbidden position remains.  By the definition of $B$, the corresponding centre vertex has a residual edge to an available cycle in each graph; because the position is non-forbidden, neither leaf has been used before.  Repeating this argument yields $t$ distinct leaves on each side.

Map the centre cycles using the rotation selected above.  For each chosen pair of leaves, choose the orientation and rotation of the second leaf so that the endpoint of the selected residual edge in $G_1$ maps to the corresponding endpoint in $G_2$. These maps preserve all factor edges and the $t$ selected centre-to-leaf residual edges. Since the chosen cycles are vertex-disjoint, the maps are compatible.
\end{proof}

\begin{lemma}[Iterating the stars]\label{lem:iterate-stars}
Under the deterministic assumptions above, for every fixed $\eta>0$ and positive integer $t$, and all sufficiently large $n$,
\begin{equation}\label{eq:star-det-bound}
M(G_1,G_2)
\ge
n+\frac{t}{t+1}(1-\eta)q-t.
\end{equation}
\end{lemma}

\begin{proof}
Start with all $q$ factor cycles available.  As long as at least $\eta q$ cycles remain, apply \cref{lem:common-star}, align one centre and $t$ leaves, and remove those $t+1$ cycles from the available set in each graph.  If $a$ stages are completed and $r=q-a(t+1)$ cycles remain, then at stopping $r<\eta q$, so
\[
a=\frac{q-r}{t+1}
>\frac{(1-\eta)q}{t+1}-1.
\]
The maps constructed at different stages act on disjoint vertex sets. Pair and align the remaining factor cycles arbitrarily. The resulting permutation preserves all $qk=n$ factor edges and gains $t$ residual edges from each completed star.  Hence
\[
M(G_1,G_2)
\ge n+ta
>
n+\frac{t}{t+1}(1-\eta)q-t,
\]
which proves~\eqref{eq:star-det-bound}.
\end{proof}

\subsection{The divisible case}\label{subsec:cycle-factor-divisible}

\begin{proof}[Proof of \cref{thm:cycle-factor-overlap}]
Fix $\eta>0$ and a positive integer $t$.  By \cref{thm:cycle-factor-existence} and~\eqref{eq:k-regime}, each graph contains a $k$-cycle factor with probability $1-o(1)$.  By \cref{prop:all-cycle-factors-good}, with probability $1-o(1)$ Property A holds and every such factor satisfies Property B$(\eta)$. Thus, with probability $1-o(1)$, each graph contains a factor to which \cref{lem:iterate-stars} applies.

Since
\[
q=\frac nk=(a_d+o(1))\frac n{\log n},
\]
\cref{lem:iterate-stars} yields
\[
M(G_1,G_2)
\ge
n+\left(\frac{t}{t+1}(1-\eta)a_d-o(1)\right)
\frac n{\log n}.
\]
Equivalently,
\[
M(G_1,G_2)
\ge
n+\left(\frac{t}{t+1}(1-\eta)-o(1)\right)\frac nk.
\]
Given $\varepsilon>0$, choose first $\eta>0$ small enough and then $t$ large enough that
\[
\frac{t}{t+1}(1-\eta)>1-\frac{\varepsilon}{2}.
\]
These constants are fixed before $n\to\infty$, and \cref{lem:common-star} eventually applies because $k\to\infty$.  The theorem follows.
\end{proof}

The star construction gains asymptotically one residual common edge per factor cycle.  This is why stars, rather than pairwise cycle matchings, recover the full coefficient $a_d$: grouping cycles only in pairs would yield roughly one extra edge for every two factor cycles and lose a factor two in the leading constant.

\subsection{Removing the divisibility condition}\label{subsec:remove-divisibility}

The remaining issue is that the cycle-factor theorem is applied to a size $N$ that is divisible by $k$, whereas the original graph size $n$ is arbitrary, subject to $dn$ being even.  The following coupling is proved in \cref{app:cycle-factor-technical}.

\begin{lemma}[Size-reduction coupling]\label{lem:size-reduction-coupling}
Fix $d$. Let $N<n$ be admissible and suppose $r:=n-N=O_d(\log n)$.  There is a coupling of two independent graphs $G_1,G_2\sim\mathcal U_{n,d}$ and two independent graphs $H_1,H_2\sim\mathcal U_{N,d}$ such that, with probability $1-o(1)$,
\begin{equation}\label{eq:size-reduction-overlap}
M(G_1,G_2)\ge M(H_1,H_2)-dr.
\end{equation}
\end{lemma}

We can now prove the lower bound stated in the introduction.

\begin{proof}[Proof of \cref{thm:second-order-lower}]
Fix $\varepsilon>0$ and choose $\varepsilon_0>0$ so small that
\[
(1-\varepsilon_0)a_d>a_d-\frac{\varepsilon}{3}.
\]
Set
\[
k=\left\lceil\frac{\log n}{a_d}+\sqrt{\log n}\right\rceil.
\]
If $d$ is odd, let $L=\operatorname{lcm}(k,2)$; if $d$ is even, let $L=k$.  Let $N$ be the largest multiple of $L$ not exceeding $n$.  Then $dN$ is even, $k\mid N$, and
\[
0\le n-N<L\le2k=O_d(\log n).
\]

Since $N=n-O_d(\log n)$, we have
\[
k=\frac{\log N}{a_d}+o(\log N).
\]
Moreover, as $N\le n$ and $a_d>0$,
our choice of $k$ ensures that
\[
k\ge
\frac{\log n}{a_d}+\sqrt{\log n}
\ge
\frac{\log N}{a_d}+\sqrt{\log N}.
\]
Thus, all the assumptions of \cref{thm:cycle-factor-overlap} are satisfied, and applying that theorem with $\varepsilon_0$ gives
\[
M(H_1,H_2)
\ge
N+(1-\varepsilon_0)\frac Nk
=
N+\left((1-\varepsilon_0)a_d-o(1)\right)\frac N{\log N}
\]
with high probability.

Apply \cref{lem:size-reduction-coupling}. Since $r=n-N=O_d(\log n)$,
\begin{align*}
M(G_1,G_2)
&\ge
N+\left((1-\varepsilon_0)a_d-o(1)\right)\frac N{\log N}-d(n-N)\\
&=
n+\left((1-\varepsilon_0)a_d-o(1)\right)\frac n{\log n},
\end{align*}
because $n-N=o(n/\log n)$ and $N/\log N=(1+o(1))n/\log n$. By the choice of $\varepsilon_0$, the last expression is at least
\[
n+(a_d-\varepsilon)\frac n{\log n}
\]
for all sufficiently large $n$, on an event of probability tending to one.  This proves the theorem.
\end{proof}

The contrast with \cref{sec:hamiltonian-alignment} can now be stated cleanly.  For a fixed realised cycle factor, the residual pairing is uniform conditional on the factor, and the exponentially small probability of violating Property B can be summed over the entire near-threshold family because that family has only $\exp(o(n))$ members in expectation.  The cycle factor may therefore be selected after observing the graph. In the Hamilton-cycle argument, we instead work with a distinguished planted cycle and use contiguity to transfer the resulting graph event.
\section{Discussion and open questions}\label{sec:discussion}

For every fixed $d\ge3$, we have shown that the maximum overlap between two independent uniform $d$-regular graphs satisfies
\[
(a_d-o(1))\frac{n}{\log n}
\le M_n - n\le
(b_d+o(1))\frac{n}{\log n}
\]
with high probability, where
\[
a_d=\log(d-1)+\frac{d-2}{2}\log\!\left(\frac{d-2}{d}\right),
\qquad
b_d=d\log(d-1)-(d-2)\log(d+1).
\]
The two bounds arise from rather different arguments. The lower bound is produced by spanning cycle factors with cycle length of order $\log n$, together with a deterministic star construction that preserves asymptotically one additional residual edge per factor cycle. The upper bound instead comes from a first-moment count of common cores with minimum degree two, and is dominated by non-spanning witnesses. Understanding the gap between these two mechanisms is the main question left open by our results.

The Hamilton-cycle argument of \cref{sec:hamiltonian-alignment} illustrates a different, complementary mechanism. There, one starts with a distinguished spanning cycle, analyses the residual graph using a planted model, and transfers the resulting event via contiguity.  For logarithmic cycle factors, the family of near-threshold factors is sufficiently small that the residual properties needed for the alignment can instead be established uniformly over all realised factors. This distinction between a planted distinguished structure and uniform control over an entire family of structures may be useful beyond the present problem.

\paragraph{Toward a sharp second-order constant.}\label{subsec:sharp-constant}
It is natural to ask whether there exists a deterministic constant $\gamma_d$ such that
\[
\frac{(M_n-n)\log n}{n}
\xrightarrow{\Pp}
\gamma_d.
\]
On the lower-bound side, the current construction is naturally limited by the number of near-threshold factor cycles. If $q=(a_d+o(1))\frac{n}{\log n}$ is the number of cycles in the factor, the quotient structure produced by the star construction is a forest and therefore contains at most $q-1$ residual edges.  In this sense, the coefficient $a_d$ is the natural endpoint of the present cycle-factor method.  Improving it would require a common residual structure denser than a forest at the quotient level, or a different spanning decomposition.

The upper bound points to a different obstruction. In \cref{sec:upper}, the first moment for common cores is maximised by witnesses occupying asymptotically $\frac{dn}{d+1}$ vertices. If one artificially restricted the count to spanning witnesses with minimum degree two, the same calculation would give the smaller coefficient $2a_d$; see~\eqref{remark:2ad}. In particular, the present union bound is limited by the entropy of non-spanning cores.

\paragraph{Beyond fixed-degree regular graphs.}\label{subsec:extensions}
The present paper deliberately restricts its attention to fixed $d$. This assumption enters at several points: the simplicity probability of the configuration model remains bounded away from zero, the Hamiltonicity and cycle-factor inputs are stated for fixed degree, and all pairing estimates are used with constants that depend on $d$. It would therefore be interesting to investigate the regime $d=d(n)\to\infty$. 
Hamiltonicity is known to hold with high probability also for non-constant degrees in the range $3 \le d \le cn$, for a sufficiently small constant $c>0$ (see \cite{cooper2002hamiltonicity}). The formulas above suggest that the second-order constants should grow like $\log d$, but the arguments given here are not uniform in $d$ and do not identify the range of growth for which the scale $n/\log n$ remains valid without modification.

A second natural direction is to replace regular graphs by configuration models with non-constant degree sequences. There are broad classes of such models, in particular under suitable assumptions with a minimum degree of at least three, for which Hamiltonicity holds with high probability \cite{cooper2010hamilton}. This suggests that a Hamiltonian baseline may survive well beyond the regular case. The main difficulty is that the two ingredients used for the stronger lower bound become degree-dependent: one needs an analogue of the near-threshold logarithmic cycle-factor theorem, and after removing such a factor, the residual half-edge distribution is no longer regular. The constants governing both the factor threshold and the residual density should then depend on the underlying degree distribution. Developing a version of the present argument in that setting would provide a natural test of how much of the phenomenon is genuinely tied to regularity.

It is also useful to contrast the regular model with sparse Erd\H{o}s--R\'enyi graphs. The methods of Sections~\ref{sec:hamiltonian-alignment} and~\ref{sec:lower} do not transfer directly to $G(n,c/n)$. In the regular model, the hard degree constraint provides a spanning Hamiltonian baseline and, after removing a spanning $2$-factor, leaves a residual graph with a particularly tractable pairing structure. In sparse Erd\H{o}s--R\'enyi graphs, low-degree and isolated vertices prevent such a spanning baseline in the same parameter range, and the residual graph obtained after exposing a spanning structure does not have the same regular-pairing form.  

\paragraph{Algorithmic questions.}\label{subsec:algorithms}

The arguments in this paper are primarily information-theoretic. Once suitable short cycle factors satisfying the required residual properties are given, the deterministic star alignment from \cref{sec:lower} is constructive: one repeatedly chooses a centre cycle with many residual attachments, tests the $k=O(\log n)$ possible relative rotations of a pair of centres, and selects a fixed number of distinct leaves. The non-constructive step in the present proof is the production of a near-threshold spanning $k$-cycle factor in a typical random regular graph.

This suggests a natural algorithmic version of the problem: find, with high probability and in polynomial time, an alignment achieving a surplus of order $n/\log n$, ideally with an explicit constant. In this regime, multiplicative approximation of the total overlap is not the most informative objective, since the Hamiltonian baseline already contributes $n$ common edges while the optimal improvement is only of order $n/\log n$. The appropriate algorithmic target is therefore the additive surplus above the baseline.

\section*{Acknowledgements}
The authors are partially funded by the MUR (Italy) Department of Excellence 2023–2027.

\section*{Conflict of interest}
None to declare.

\section*{Data availability}
No datasets were generated or analysed in this theoretical study.

\section*{AI usage statement}
This project fully aligns with the European Commission's \emph{ERA Living Guidelines on the Responsible Use of Generative AI in Research} (October 2026). The authors disclose the use of ChatGPT6 (OpenAI) and Grammarly during the preparation of this manuscript, for language editing, to improve exposition, and to organise mathematical arguments.

In particular, ChatGPT6 assisted in identifying and resolving an issue in the proof of the upper bound in \cref{thm:overlap-d2}, arising from difficulties in applying a subset-sum estimate of Brito et al. It also suggested two improvements to the leading order constant of Theorems \ref{thm:second-order-upper} and \ref{thm:second-order-lower}: a slight refinement of the upper bound, and a more substantial improvement of the lower-bound constant by a factor of two. The latter was achieved by replacing an earlier contiguity-based approach with the actual argument described in \cref{sec:lower}.

All such mathematical arguments were subsequently developed, independently checked and verified by the authors, who take full responsibility for the statements, proofs, and conclusions of the manuscript.

\bibliographystyle{amsplain}
\bibliography{references}

\newpage
\appendix
\section{Common subset sums of Ewens partitions}\label{app:ewens-common-sums}

In this section, we prove the following lemma.

\begin{lemma}[Common subset sums]\label{lem:degree-two-subset-sums}
Let $G_1,G_2$ be independent samples from $\mathcal U_{n,2}$, and let $\mathcal L_i$ be the set of all subset sums of the cycle lengths of $G_i$, including $0$ and $n$. Then
\[
\bigl|(\mathcal L_1\cap\mathcal L_2)\setminus\{0,n\}\bigr|=O_{\Pp}(1).
\]
\end{lemma}
\noindent Combining this result with \eqref{eq:bZ} yields a sharper upper bound in the proof of \cref{thm:overlap-d2}.

Recall that the cycle counts of a uniform simple $2$-regular graph follow the Ewens sampling distribution with parameter $\theta=1/2$, conditioned on having no cycles of lengths $1$ or $2$. For background on the Ewens sampling formula and the Feller coupling, see \cite{tavare2021magical}.

Our approach builds on the study of subset sums of Ewens cycle lengths; see \cite{brito2018}. We separate small and macroscopic subset sums. For small sums, the Feller coupling reduces the problem to the intersection of two independent Poisson sumsets, which we show to be almost surely finite. For macroscopic sums, a weighted generating-function argument shows that two independent Ewens partitions have no common subset sums with high probability. Together, these estimates yield the required tightness.

For a sequence $x=(x_\ell)_{\ell\geq1}$ of non-negative integers, define its sumset by
\[
\mathcal L(x):=\left\{\sum_{\ell\geq1}\ell a_\ell:\;0\leq a_\ell\leq x_\ell,\ \text{with finite support}\right\}.
\]
For a cycle-count vector $C^{(n)}$, the set $\mathcal L(C^{(n)})$ consists precisely of the integers that can be obtained by summing the lengths of a collection of cycles. Let $Y=(Y_\ell)_{\ell\geq1}$ have independent coordinates with
\[
Y_\ell\sim\operatorname{Poi}\!\left(\frac{1}{2\ell}\right).
\]
We first establish a property of the infinite sumsets that will be useful under the Feller coupling.

\begin{lemma}\label{lem:poisson-sumset-intersection}
If $Y^{(1)},Y^{(2)}$ are independent copies of $Y$, then
\[
T:=\bigl|\mathcal L(Y^{(1)})\cap\mathcal L(Y^{(2)})\bigr|<\infty
\qquad\text{almost surely}.
\]
\end{lemma}

\begin{proof}

Fix $a=11/20$. This choice satisfies
\[
\frac12<a<
\min\left\{
\frac{1}{2\log 2},
\frac{1}{4\log(3/2)}
\right\},
\]
as required by the summability estimates in the two parts of the argument. For $m\geq1$ put
\[
N_m=\sum_{\ell=1}^mY_\ell,
\qquad
B_m=[z^m]\prod_{\ell=1}^m(1+z^\ell)^{Y_\ell}.
\]
Thus $B_m$ counts the selections of available summands whose total is $m$, and $m\in\mathcal L(Y)$ implies $B_m\geq1$. Writing $H_m=\sum_{\ell=1}^m\ell^{-1}$, independence and the Poisson generating function give
\begin{align*}
\E{2^{-N_m}B_m}
&=[z^m]\prod_{\ell=1}^m
   \E{\left(\frac{1+z^\ell}{2}\right)^{Y_\ell}}\\
&=e^{-H_m/4}[z^m]\exp\!\left(\frac14\sum_{\ell=1}^m\frac{z^\ell}{\ell}\right)\\
&=e^{-H_m/4}[z^m](1-z)^{-1/4}
 =O(m^{-1}).
\end{align*}
The third equality holds because terms of degree greater than $m$ do not affect the coefficient of $z^m$. On $\{N_m\leq a\log m\}$,
$2^{-N_m}B_m\geq m^{-a\log2}$ whenever $m\in\mathcal L(Y)$. Hence
\begin{equation}\label{eq:poisson-quenched-bound}
\Pp\bigl(m\in\mathcal L(Y),\ N_m\leq a\log m\bigr)
 \leq C m^{-1+a\log2}.
\end{equation}
By independence, the probability that $m$ belongs to both sumsets and that both corresponding bounds on $N_m$ hold is at most
\[
C^2m^{-2+2a\log2}.
\]
The sum of these probabilities over $m\geq2$ is finite, since $2a\log2<1$.

Finally, $N_m$ has the same finite-dimensional distributions as a unit-rate Poisson process sampled at time $H_m/2$. Consequently $N_m/\log m\to1/2$ almost surely, and $N_m\leq a\log m$ eventually almost surely. The same holds for the independent copy. Borel--Cantelli applied to the summable probabilities above proves that only finitely many positive integers belong to both sumsets almost surely. Including $0$ does not affect finiteness.
\end{proof}

Now, let $C^{(n)}=(C_1^{(n)},\ldots,C_n^{(n)})$ have the Ewens sampling distribution $\operatorname{ESF}(1/2,n)$, and set $\mathcal L_n=\mathcal L(C^{(n)})$ and $\kappa_n=\sum_{\ell=1}^nC_\ell^{(n)}$. For $t>0$ and $m\geq 0$, define
\[
h_m(t):=[z^m](1-z)^{-t/2}
=\frac{\Gamma(m+t/2)}{\Gamma(t/2)\Gamma(m+1)}.
\]
The Ewens sampling formula gives
\[
\Pp(C^{(n)}=c)
=\frac{1}{h_n(1)}\prod_{\ell=1}^n\frac{(1/(2\ell))^{c_\ell}}{c_\ell!},
\qquad \sum_{\ell=1}^n\ell c_\ell=n.
\]
For $1\leq k\leq n-1$, let $A_k$ denote the number of collections of cycles whose total length is $k$, so that $\mathbf{1}_{\{k\in\mathcal L_n\}}\leq A_k$. Marking the cycles selected by such a subset and those in its complement separately gives the exact identity
\begin{equation}\label{eq:ewens-weighted-subsets}
\E{A_k t^{\kappa_n}}=\frac{h_k(t)h_{n-k}(t)}{h_n(1)}.
\end{equation}

\begin{lemma}\label{lem:ewens-middle-subsets}
If $\mathcal L_n^{(1)},\mathcal L_n^{(2)}$ arise from two independent $\operatorname{ESF}(1/2,n)$ permutations, then
\[
\Pp\left(
\bigl(\mathcal L_n^{(1)}\cap\mathcal L_n^{(2)}\bigr)
\cap[\lfloor\sqrt n\rfloor,n-\lfloor\sqrt n\rfloor]
\neq\varnothing
\right)=o(1).
\]
\end{lemma}

\begin{proof}
Again set $a=11/20$, and let $\mathcal E_n=\{\kappa_n\leq a\log n\}$. The classical independent-Bernoulli representation of the total number of cycles under Ewens' distribution implies
\[
\E{\kappa_n}=\frac12\log n+O(1),
\qquad\operatorname{Var}(\kappa_n)=O(\log n),
\]
so $\Pp(\mathcal E_n)=1-o(1)$ by Chebyshev's inequality.

Take $t=2/3$. On $\mathcal E_n$ we have
$\mathbf1_{\mathcal E_n}\leq t^{-a\log n}t^{\kappa_n}$.
For $1\leq k\leq n/2$, using \eqref{eq:ewens-weighted-subsets} and the coefficient estimates
$h_m(2/3)=O((m+1)^{-2/3})$ and
$h_n(1)^{-1}=O(n^{1/2})$, we obtain
\begin{align}
\Pp(k\in\mathcal L_n,\mathcal E_n)
&\leq t^{-a\log n}\E{A_k t^{\kappa_n}}\notag\\
&\leq Cn^{a\log(3/2)-1/6}k^{-2/3}.
\label{eq:ewens-macro-prob}
\end{align}
Here we used $n-k\geq n/2$, and the constant is independent of $k$ and $n$.

Let $R=\lfloor\sqrt n\rfloor$. By the independence of the two Ewens partitions and a union bound,
\begin{align*}
&\Pp\Bigl(\exists k\in[R,\lfloor n/2\rfloor]:
 k\in\mathcal L_n^{(1)}\cap\mathcal L_n^{(2)},\;
 \mathcal E_n^{(1)}\cap\mathcal E_n^{(2)}\Bigr)\\
&\quad\leq Cn^{2a\log(3/2)-1/3}
\sum_{k=R}^{\lfloor n/2\rfloor}k^{-4/3}\\
&\quad\leq Cn^{2a\log(3/2)-1/2}=o(1),
\end{align*}
because $2a\log(3/2)<1/2$. Finally, $k\in\mathcal L_n$ if and only if $n-k\in\mathcal L_n$, since taking the complementary set of cycles preserves the Ewens partition. The estimate, therefore, extends to $[R,n-R]$. Adding $\Pp((\mathcal E_n^{(1)})^c\cup(\mathcal E_n^{(2)})^c)=o(1)$ proves the lemma.
\end{proof}

We are now ready to conclude the proof of the main lemma.

\begin{proof}[Proof of \cref{lem:degree-two-subset-sums}]
First consider two independent Ewens partitions $C^{(n,1)},C^{(n,2)}$ with parameter $1/2$; write $\mathcal L_n^{(i)}=\mathcal L(C^{(n,i)})$. Let $R=\lfloor\sqrt n\rfloor$.

By the Feller coupling (see \cite[Section~3.1, Lemma~10]{brito2018}), each $C^{(n,i)}$ can be coupled to a sequence $Y^{(i)}$ of independent $\operatorname{Poi}(1/(2\ell))$ variables in such a way that
\[
C_\ell^{(n,i)}\leq Y_\ell^{(i)}+\mathbf1_{\{J_n^{(i)}=\ell\}},
\qquad 1\leq\ell\leq n,
\]
where the same reference gives, for $\ell<n$,
\[
\Pp(J_n^{(i)}=\ell)\leq\frac{1}{2(n-\ell)}.
\]
Consequently,
\[
\Pp(J_n^{(1)}\leq R\text{ or }J_n^{(2)}\leq R)
\leq\sum_{i=1}^2\sum_{\ell=1}^{R}\frac{1}{2(n-\ell)}
=O(n^{-1/2}).
\]
On the complementary event, for $i=1,2$,
\[
\mathcal L_n^{(i)}\cap[0,R]
\subseteq\mathcal L(Y^{(i)})\cap[0,R],
\]
since every subset sum at most $R$ uses only cycles of length at most $R$.

By \cref{lem:ewens-middle-subsets}, with probability $1-o(1)$ there are no common subset sums in $[R,n-R]$. By symmetry around $n/2$, the number of remaining non-trivial common subset sums is at most
\[
2\bigl|\mathcal L(Y^{(1)})\cap\mathcal L(Y^{(2)})\bigr|.
\]
The quantity on the right is almost surely finite by \cref{lem:poisson-sumset-intersection}, and its law does not depend on $n$. This shows that
\[
\bigl|(\mathcal L_n^{(1)}\cap\mathcal L_n^{(2)})\setminus\{0,n\}\bigr|
=O_{\Pp}(1)
\]
for two independent Ewens partitions.

Finally, the cycle counts of a uniform simple $2$-regular graph have exactly the $\operatorname{ESF}(1/2,n)$ distribution conditioned on $C_1^{(n)}=C_2^{(n)}=0$. The probability of this event tends to $e^{-3/4}>0$ by the convergence of the first two Ewens cycle counts to independent Poisson variables with means $1/2$ and $1/4$. Conditioning two independent Ewens partitions on these two events thus preserves the $O_{\Pp}(1)$ bound (at a bounded multiplicative cost in probabilities), and gives the statement for $G_1,G_2\sim\mathcal U_{n,2}$.
\end{proof}

\section{Technical estimates for the upper bound}\label{app:upper-bound-count}

This appendix proves \cref{prop:core-first-moment}.  We keep the notation of \cref{sec:upper}. Thus, $v$ is the number of vertices of the common core, $s$ is its excess, and
\[
    m=v+s.
\]
All estimates below are for fixed $ d\ge 3$.  We work first in the configuration model; conditioning on simplicity is performed only in the proof of \cref{thm:second-order}.

\subsection{Configuration-model enumeration}

\begin{lemma}[Counting common cores]\label{lem:app-core-count}
For every feasible $v$ and $s$,
\begin{align}
\Ee Z_{v,s}
\le{}&
\binom nv^2 v!\,(2m-1)!!
\left(
    \frac{(dn-2m-1)!!}{(dn-1)!!}
\right)^2
\nonumber\\
&\times
\sum_{\substack{a_1,\ldots,a_v\in\{2,\ldots,d\}\\
                  \sum_i a_i=2m}}
\prod_{i=1}^v\frac{(d)_{a_i}^{\,2}}{a_i!}.
\label{eq:app-core-count}
\end{align}
Here $(d)_r=d(d-1)\cdots(d-r+1)$.
\end{lemma}

\begin{proof}
Choose the two $v$-vertex sets and the bijection between them in
\[
    \binom nv^2v!
\]
ways. Fix a degree sequence $a_1,\ldots,a_v\in\{2,\ldots,d\}$ with sum $2m$.

The number of simple labelled graphs with this degree sequence is at most
\[
    \frac{(2m-1)!!}{\prod_{i=1}^v a_i!}.
\]
Indeed, attach $a_i$ labelled auxiliary half-edges to vertex $i$ and pair the resulting $2m$ half-edges.  Every simple graph with the prescribed degree sequence is generated by exactly $\prod_i a_i!$ such pairings, whereas pairings that produce loops or multiple edges only increase the upper bound.

Now fix one such simple graph $H$.  At vertex $i$, the $a_i$ incidences of $H$ can be assigned injectively to the $d$ configuration half-edges in $(d)_{a_i}$ ways.  Once all incidence assignments are fixed, the $m$ edges of $H$ prescribe $m$ disjoint half-edge pairs.  A uniform pairing of the $dn$ half-edges contains these prescribed pairs with probability
\[
    \frac{(dn-2m-1)!!}{(dn-1)!!}.
\]
Consequently,
\[
    \Pp(H\subseteq G^{\mathrm{CM}}_{n,d})
    \le
    \left(\prod_{i=1}^v(d)_{a_i}\right)
    \frac{(dn-2m-1)!!}{(dn-1)!!}.
\]
The two configuration models are independent, and the image of $H$ under the chosen bijection has the same degree sequence.  Squaring the containment probability, multiplying by the graph count, and summing over all admissible degree sequences gives \eqref{eq:app-core-count}.
\end{proof}

\subsection{The entropy of the degree sequence}

The minimum-degree condition makes degree two the natural baseline.  Write
\[
    a_i=2+b_i,
    \qquad
    0\le b_i\le d-2.
\]
Since $\sum_i a_i=2v+2s$, the excess is encoded by
\[
    \sum_{i=1}^v b_i=2s.
\]

\begin{lemma}[Degree-sequence entropy]\label{lem:app-degree-entropy}
Set
\[
    \Lambda_d:=\frac{d^2(d-1)^2}{2}.
\]
Uniformly for $1\le s\le n/2$ and $1\le v\le n$,
\begin{equation}\label{eq:app-degree-sum}
\sum_{\substack{a_1,\ldots,a_v\in\{2,\ldots,d\}\\
                  \sum_i a_i=2v+2s}}
\prod_{i=1}^v\frac{(d)_{a_i}^{\,2}}{a_i!}
\le
\Lambda_d^v
\exp\!\left\{
    O_d\!\left(s\log\frac ns+s\right)
\right\}.
\end{equation}
\end{lemma}

\begin{proof}
For $0\le b\le d-2$, put
\[
    w_b:=\frac{(d)_{2+b}^2}{(2+b)!}.
\]
Then $w_0=\Lambda_d$. Since $d$ is fixed, there is a constant $C_d\ge1$ such that
\[
    w_b\le \Lambda_d C_d^b,
    \qquad 0\le b\le d-2.
\]
Hence every admissible vector $(b_1,\ldots,b_v)$ with sum $2s$ contributes at most
\[
    \Lambda_d^v C_d^{2s}.
\]

Ignoring the upper bounds $b_i\le d-2$ can only increase the number of vectors. The number of nonnegative integer solutions to $\sum_i b_i=2s$ is
\[
    \binom{v+2s-1}{2s}
    \le
    \left(\frac{e(v+2s)}{2s}\right)^{2s}
    \le
    \exp\!\left\{
        O\!\left(s\log\frac ns+s\right)
    \right\},
\]
where we used $v\le n$ and $s\le n/2$. Multiplying the two estimates proves \eqref{eq:app-degree-sum}.
\end{proof}

\subsection{The cost of positive excess}

Define the excess-zero comparison term
\begin{equation}\label{eq:app-Av}
    A_v
    :=
    \binom nv^2v!\,(2v-1)!!
    \left(
        \frac{(dn-2v-1)!!}{(dn-1)!!}
    \right)^2
    \Lambda_d^v.
\end{equation}

\begin{lemma}[One factor $n^{-1}$ per unit of excess]\label{lem:app-excess-cost}
If $s=O(n/\log n)$, then uniformly over feasible $v$,
\begin{equation}\label{eq:app-compare-s0}
    \Ee Z_{v,s}
    \le
    A_v
    \exp\!\left\{
        -s\log n
        +O_d\!\left(s\log\frac ns+s\right)
    \right\}.
\end{equation}
\end{lemma}

\begin{proof}
By Lemmas~\ref{lem:app-core-count} and~\ref{lem:app-degree-entropy}, it remains to compare the pairing factor for $m=v+s$ with the corresponding factor for $m=v$. Their ratio is
\begin{align*}
&\frac{(2(v+s)-1)!!}{(2v-1)!!}
\left(
    \frac{(dn-2(v+s)-1)!!}{(dn-2v-1)!!}
\right)^2\\
&\hspace{2cm}=
\prod_{j=0}^{s-1}
\frac{2v+2j+1}{(dn-2v-2j-1)^2}.
\end{align*}
For $v\le n$ and $s=O(n/\log n)$, the numerator of every factor is $O(n)$, while
\[
    dn-2v-2j-1
    \ge
    (d-2)n-O(n/\log n)
    \ge c_dn
\]
for some $c_d>0$ and all sufficiently large $n$. Hence, each factor is at most $C_d/n$, uniformly in $v$ and $j$.  The whole ratio is therefore at most
\[
    n^{-s}e^{O_d(s)}.
\]
Combining this estimate with \eqref{eq:app-degree-sum} proves \eqref{eq:app-compare-s0}.
\end{proof}

\subsection{Entropy of the excess-zero core}

\begin{lemma}[Entropy of the core size]\label{lem:app-core-entropy}
Let $x=v/n$.  Uniformly for $1\le v\le n$,
\begin{equation}\label{eq:app-A-f}
    \log A_v
    =n f_d(x)+O_d(\log n),
\end{equation}
where $f_d$ is defined in \eqref{eq:f-def}.
\end{lemma}

\begin{proof}
Use
\[
    (2r-1)!!=\frac{(2r)!}{2^r r!}
\]
and Stirling's formula in the form
\[
    \log r!=r\log r-r+O(\log(r+1)).
\]
Uniformly for $1\le v\le n$, with $x=v/n$,
\begin{align}
\log\!\left[\binom nv^2v!\right]
&=
2n[-x\log x-(1-x)\log(1-x)]
+xn\log(xn)-xn
+O(\log n),
\label{eq:app-choose}\\
\log(2v-1)!!
&=
 xn\log(2xn)-xn+O(\log n).
\label{eq:app-df-v}
\end{align}
For the ratio of double factorials, write
\begin{align*}
\log\frac{(dn-2v-1)!!}{(dn-1)!!}
&=-\sum_{j=0}^{v-1}\log(dn-2j-1)\\
&=-v\log n
  -\sum_{j=0}^{v-1}
     \log\!\left(d-\frac{2j+1}{n}\right).
\end{align*}
Since $d\ge3$, the function $y\mapsto\log(d-2y)$ and its derivative are uniformly bounded on $[0,1]$. A Riemann-sum estimate therefore gives
\begin{equation}\label{eq:app-ratio-df}
\log\frac{(dn-2v-1)!!}{(dn-1)!!}
=
-xn\log n
-n\int_0^x\log(d-2y)\,\dd y
+O_d(1).
\end{equation}
The integral is
\begin{equation}\label{eq:app-integral}
\int_0^x\log(d-2y)\,\dd y
=
\frac12\left[
    d\log d-(d-2x)\log(d-2x)-2x
\right].
\end{equation}
Finally,
\[
    v\log\Lambda_d
    =xn\log\frac{d^2(d-1)^2}{2}.
\]
Substituting \eqref{eq:app-choose}--\eqref{eq:app-integral} into \eqref{eq:app-Av}, all $n\log n$ terms cancel and the remaining terms simplify to
\[
    n\left[
        -2(1-x)\log(1-x)
        +2x\log[d(d-1)]
        -d\log d
        +(d-2x)\log(d-2x)
    \right]
    +O_d(\log n),
\]
which is exactly \eqref{eq:app-A-f}.
\end{proof}

\begin{proof}[Proof of \cref{prop:core-first-moment}]
Combining \cref{lem:app-excess-cost} with \cref{lem:app-core-entropy} gives
\[
    \log \Ee Z_{v,s}
    \le
    n f_d(v/n)-s\log n
    +O_d\!\left(s\log\frac ns+s+\log n\right).
\]
In the regime $1\le s=O(n/\log n)$, the additive $O_d(s)$ term is absorbed by
\[
    O_d\!\left(s\log\frac ns\right),
\]
for all sufficiently large $n$, because $n/s\to\infty$. This proves \eqref{eq:core-first-moment}.
\end{proof}

\section{Technical estimates for Hamilton-cycle alignments}
\label{app:hamiltonian-technical}

This appendix contains the auxiliary matching estimates and second-moment calculations used in \cref{sec:hamiltonian-alignment}.  We also give the formal proof of \cref{thm:typical-hamiltonian-alignment} and the upper bound needed for \cref{prop:cyclic-matching-overlap} and \cref{thm:cyclic-matching-unions}.

\subsection{Fixed subgraphs, extracted matchings, and edge-disjointness}
\label{app:ham-basic-matching}

We first collect three elementary facts about uniform regular graphs.

\begin{lemma}[Fixed-subgraph containment]
\label{lem:fixed-subgraph}
Fix $q\ge1$ and a positive integer $j$. Let $F\subseteq\binom{[n]}2$ be a deterministic set of $j$ distinct edges, and write $r_v=d_F(v)$. If $r_v\le q$ for every vertex, then for $R\sim\mathcal U_{n,q}$,
\[
\Pp(F\subseteq E(R))
=(1+o(1))\frac{\prod_{v=1}^n(q)_{r_v}}{(qn)^j},
\]
uniformly over the choice of $F$. If some $r_v>q$, the containment probability is zero. In particular, if $F$ is a matching of $j$ edges, then
\[
\Pp(F\subseteq E(R))
=(1+o(1))\left(\frac qn\right)^j.
\]
\end{lemma}

\begin{proof}
In the configuration model, the half-edges corresponding to the incidences of $F$ can be assigned in $\prod_v(q)_{r_v}$ ways. For each assignment, the probability that the resulting $j$ prescribed pairs occur is
\[
\prod_{t=0}^{j-1}\frac1{qn-2t-1}.
\]
Since $j$ and $q$ are fixed,
\[
\Pp(F\subseteq E(G^{\mathrm{CM}}_{n,q}))
=(1+o(1))\frac{\prod_v(q)_{r_v}}{(qn)^j}.
\]
Conditioning on simplicity preserves this asymptotic estimate for fixed prescribed subgraphs; see~\cite{Wormald1999}.  The matching case follows because exactly $2j$ vertices have degree one in $F$.
\end{proof}

\begin{lemma}[Uniform matching extraction]
\label{lem:matching-extraction}
Fix $q\ge1$ and let $R\sim\mathcal U_{n,q}$. Set $a=\lfloor n/4\rfloor$. There exists a random matching $J\subseteq E(R)$ such that $|J|=a$ and $J$ is uniformly distributed among all matchings of size $a$ on $[n]$.
\end{lemma}

\begin{proof}
Every $q$-regular graph contains a matching of size at least $n/4$. Indeed, a greedy maximal matching selects one edge at a time and removes at most $2q-1$ edges at each step. Since the graph has $qn/2$ edges, the resulting matching has size at least
\[
\frac{qn/2}{2q-1}\ge\frac n4.
\]
Conditionally on $R$, choose uniformly a matching $J\subseteq E(R)$ of size $a$.  The law of $R$ and the selection rule are invariant under every permutation of the vertices. Since the symmetric group acts transitively on matchings of a fixed size, the marginal distribution of $J$ is uniform among all size-$a$ matchings on $[n]$.
\end{proof}

\begin{lemma}[Edge-disjointness of a $2$-factor and a regular graph]
\label{lem:factor-disjointness}
Fix $q\ge1$. Let $F$ be any deterministic simple $2$-factor on $[n]$ and let $R\sim\mathcal U_{n,q}$ be independent of $F$. Then
\[
|E(F)\cap E(R)|\xrightarrow{d}\operatorname{Poisson}(q),
\]
uniformly over the choice of $F$. In particular,
\[
\Pp(E(F)\cap E(R)=\varnothing)\longrightarrow e^{-q}>0.
\]
\end{lemma}

\begin{proof}
Let $Z=|E(F)\cap E(R)|$ and fix a positive integer $j$. Since $F$ has exactly $n$ edges and maximum degree two, the number of ordered $j$-tuples of distinct, pairwise vertex-disjoint edges of $F$ is
\[
n^j+O_j(n^{j-1}),
\]
uniformly over the choice of $F$.  By \cref{lem:fixed-subgraph}, the probability that any such tuple is contained in $R$ is
\[
(1+o(1))\left(\frac qn\right)^j.
\]
There are only $O_j(n^{j-1})$ remaining ordered tuples of distinct edges of $F$, and their containment probabilities are at most $O_{q,j}(n^{-j})$, again by \cref{lem:fixed-subgraph}, with probability zero whenever a prescribed degree exceeds $q$.  Consequently,
\[
\Ee[(Z)_j]
=\left(n^j+O_j(n^{j-1})\right)(1+o(1))\left(\frac qn\right)^j+O_{q,j}(n^{-1})
\longrightarrow q^j.
\]
These are the factorial moments of a Poisson random variable with mean $q$, proving the claim.
\end{proof}

\subsection{Rotations of uniform matchings}
\label{app:ham-random-matchings}

We now prove \cref{lem:cyclic-matching-witness}.  Throughout this subsection, vertices are identified with elements of $\mathbb Z_n$.

\begin{lemma}[Containment probability for a random matching]
\label{lem:matching-containment}
Let $J$ be a uniform matching of size $a$ on $\mathbb Z_n$, and let $S$ be a fixed matching consisting of $j$ edges. Then
\[
\Prob{S\subseteq J}=p_j^{(a)}:=\frac{2^j(a)_j}{(n)_{2j}}.
\]
In particular, if $a\ge cn$ for some fixed $c>0$, and $j=o(\sqrt n)$, then, writing $\lambda_a=2a/n$,
\[
p_j^{(a)}
=\left(\frac{\lambda_a}{n}\right)^j
\exp\left(O\left(\frac{j^2}{n}\right)\right),
\]
uniformly over $cn\le a\le n/2$.
\end{lemma}

\begin{proof}
The number of matchings of size $a$ on $n$ vertices is
\[
\frac{n!}{(n-2a)!\,2^a a!}.
\]
After prescribing the $j$ edges of $S$, the remaining $a-j$ edges form an arbitrary matching on the other $n-2j$ vertices. Taking the ratio gives
\[
p_j^{(a)}=\frac{2^j(a)_j}{(n)_{2j}}.
\]
Moreover,
\[
p_j^{(a)}
=\left(\frac{\lambda_a}{n}\right)^j
\frac{\prod_{r=0}^{j-1}(1-r/a)}{\prod_{r=0}^{2j-1}(1-r/n)}.
\]
Since $a\ge cn$ and $j=o(\sqrt n)$,
\[
\log\left(\frac{p_j^{(a)}}{(\lambda_a/n)^j}\right)
=O(j^2/n),
\]
uniformly over the allowed values of $a$.
\end{proof}

For an edge $e=\{u,v\}$, define its cyclic length by
\[
\ell(e)=\min\{(v-u)\bmod n,(u-v)\bmod n\}.
\]
Thus $\ell(e)\in\{1,\ldots,n/2\}$ when $n$ is even, and rotations preserve cyclic length.

\begin{lemma}[Linearly many distinct cyclic lengths]
\label{lem:distinct-cyclic-lengths}
Fix $c\in(0,1/2]$, and let $J$ be a uniform matching of size $a$ on $\mathbb Z_n$, where
\[
cn\le a\le n/2.
\]
With probability $1-o(1)$, $J$ contains a matching $F$ satisfying
\[
|F|\ge\frac{cn}{4},
\]
such that all edges of $F$ have pairwise distinct cyclic lengths, none equal to $n/2$. Consequently, the rotations
\[
F_s:=\rho_s(F),\qquad s\in\mathbb Z_n,
\]
are pairwise edge-disjoint.
\end{lemma}

\begin{proof}
For $1\le\ell<n/2$, let $N_\ell$ be the number of edges of cyclic length $\ell$ contained in $J$. Select one edge of every represented length and denote the resulting matching by $F$. Thus
\[
|F|=m:=\sum_{\ell=1}^{\lceil n/2\rceil-1}\mathbf 1_{\{N_\ell>0\}}.
\]
For each $\ell<n/2$, there are exactly $n$ edges of length $\ell$, so \cref{lem:matching-containment} gives
\[
\Ee[N_\ell]=np_1^{(a)}=\frac{2a}{n-1}.
\]
Among the $\binom n2$ pairs of candidate edges of length $\ell$, exactly $n$ pairs share a vertex. Hence
\[
\Ee\binom{N_\ell}{2}
=\left(\binom n2-n\right)p_2^{(a)}
=\frac{2a(a-1)}{(n-1)(n-2)}.
\]
Using $\mathbf 1_{\{z>0\}}\ge z-\binom z2$ for every non-negative integer $z$,
\begin{equation}
\label{eq:expected-distinct-cyclic-lengths-app}
\Ee[m]
\ge
\left\lceil\frac n2-1\right\rceil
\left[
\frac{2a}{n-1}-\frac{2a(a-1)}{(n-1)(n-2)}
\right]
\ge\frac{cn}{2}-O(1).
\end{equation}

To obtain concentration, generate $J$ by taking a uniform permutation $\sigma$ of $\mathbb Z_n$ and pairing
\[
(\sigma_1,\sigma_2),(\sigma_3,\sigma_4),\ldots,(\sigma_{2a-1},\sigma_{2a}).
\]
Let
\[
Y_i=\Ee[m\mid \sigma_1,\ldots,\sigma_i],\qquad 0\le i\le n.
\]
After fixing $\sigma_1,\ldots,\sigma_{i-1}$, completions with $\sigma_i=x$ and $\sigma_i=y$ can be paired by transposing $x$ and $y$. Such a transposition changes at most two matching edges, hence changes $m$ by at most two. Thus $|Y_i-Y_{i-1}|\leq 2$, and the Azuma--Hoeffding inequality (see \cite{mcdiarmid1989bounded}) gives
\[
\Prob{m\le\Ee[m]-t}\le\exp\left(-\frac{t^2}{8n}\right).
\]
By \eqref{eq:expected-distinct-cyclic-lengths-app},
\[
\Prob{m<cn/4}=o(1).
\]

Finally, suppose that $\rho_s(e)=\rho_t(f)$ for $e,f\in F$. Rotations preserve cyclic length, so $\ell(e)=\ell(f)$; by construction, this forces $e=f$. Hence $\rho_{t-s}(e)=e$. A non-zero rotation can fix an unordered edge only when the edge is antipodal, i.e. has length $n/2$, which was excluded. Therefore $s=t$.
\end{proof}

\begin{proof}[Proof of \cref{lem:cyclic-matching-witness}]
Apply \cref{lem:distinct-cyclic-lengths} to $J_1$. With probability $1-o(1)$, there exists a matching $F\subseteq J_1$ with
\[
m:=|F|\ge cn/4
\]
whose edges have distinct non-antipodal cyclic lengths. Condition on any such realisation and write
\[
F_s=\rho_s(F),\qquad s\in\mathbb Z_n.
\]
For a fixed rotation $s$, call a set $S\subseteq F_s$ a $k$-witness if $|S|=k$ and $S\subseteq J_2$. Define
\[
Z_s=\sum_{\substack{S\subseteq F_s\\|S|=k}}\mathbf 1_{\{S\subseteq J_2\}},
\qquad
Z=\sum_{s\in\mathbb Z_n}Z_s.
\]
Then $Z>0$ implies $\max_s|\rho_s(J_1)\cap J_2|\ge k$.

For brevity, write
\[
p_j=p_j^{(a_2)},\qquad \lambda=\frac{2a_2}{n}\in[2c,1].
\]
By \cref{lem:matching-containment},
\[
\mu:=\Ee[Z]=n\binom mk p_k.
\]
We estimate the second moment by separating different and equal shifts.

\paragraph{Different shifts.}
Since $F_s\cap F_t=\varnothing$ whenever $s\neq t$, two witnesses from different shifts share no edges. Their union either consists of $2k$ disjoint edges, in which case its probability of being contained in $J_2$ is $p_{2k}$, or contains two incident edges, in which case the probability is zero. Hence
\[
\sum_{s\neq t}\Ee[Z_sZ_t]
\le n(n-1)\binom mk^2p_{2k},
\]
and therefore
\begin{equation}
\label{eq:different-shifts-second-moment-app}
\frac{\sum_{s\neq t}\Ee[Z_sZ_t]}{\mu^2}
\le
\left(1-\frac1n\right)\frac{p_{2k}}{p_k^2}
=1+o(1),
\end{equation}
because $k^2/n=o(1)$.

\paragraph{The same shift.}
Fix $s$. Two $k$-subsets $S,S'\subseteq F_s$ may share edges. Classifying them by $r=|S\cap S'|$ gives
\[
\Ee[Z_s^2]
=\binom mk\sum_{r=0}^k\binom kr\binom{m-k}{k-r}p_{2k-r}.
\]
For a fixed $S$, let
\[
b_r:=\frac{\binom kr\binom{m-k}{k-r}}{\binom mk},
\]
so that a union bound gives
\[
b_r\le\binom kr\left(\frac{k}{m-k+1}\right)^r.
\]
Moreover, \cref{lem:matching-containment} yields
\[
\frac{p_{2k-r}}{p_k^2}
\le
\left(\frac n\lambda\right)^r\exp(O(k^2/n)).
\]
Since $m\ge cn/4$ and $\lambda\ge2c$, for large $n$,
\[
\frac{n}{\lambda(m-k+1)}\le\frac4{c^2}.
\]
Writing $C=4/c^2$,
\begin{align}
\frac{\sum_s\Ee[Z_s^2]}{\mu^2}
&=\frac1n\sum_{r=0}^k b_r\frac{p_{2k-r}}{p_k^2}\notag\\
&\le\frac{\exp(O(k^2/n))}{n}\sum_{r=0}^k\binom kr(Ck)^r\notag\\
&=\frac{\exp(O(k^2/n))}{n}(1+Ck)^k.
\label{eq:same-shift-second-moment-app}
\end{align}
For
\[
k=\left\lfloor(1-\varepsilon)\frac{\log n}{\log\log n}\right\rfloor,
\]
we have
\[
\log\left(\frac{(1+Ck)^k}{n}\right)
=-\log n+k\log k+O(k)
=(-\varepsilon+o(1))\log n,
\]
so the right-hand side of \eqref{eq:same-shift-second-moment-app} tends to zero. Combining this with \eqref{eq:different-shifts-second-moment-app},
\[
\frac{\Ee[Z^2]}{\Ee[Z]^2}=1+o(1),
\qquad
\frac{\Var{Z}}{\Ee[Z]^2}=o(1).
\]
Chebyshev's inequality gives $\Prob{Z=0}=o(1)$. The initial distinct-length event fails with probability $o(1)$ by \cref{lem:distinct-cyclic-lengths}, completing the proof.
\end{proof}

\subsection{Typical Hamilton-cycle pairs}
\label{app:ham-typical-pairs}

We now formalise and prove \cref{thm:typical-hamiltonian-alignment}. Fix $\varepsilon\in(0,1)$ and write
\[
k_n=\left\lfloor(1-\varepsilon)\frac{\log n}{\log\log n}\right\rfloor.
\]
For a graph $G$, let $\mathcal H(G)$ denote its set of unoriented Hamilton cycles. Given graphs $G_1,G_2$ and a pair
\[
(H_1,H_2)\in\mathcal H(G_1)\times\mathcal H(G_2),
\]
call $(H_1,H_2)$ \emph{good} if, after choosing either relative orientation of the two cycles, some cyclic shift of the induced identification satisfies
\[
|\pi(E(G_1))\cap E(G_2)|\ge n+k_n.
\]
Define
\[
A_n(G_1,G_2)
:=
\frac{
\#\{(H_1,H_2)\in\mathcal H(G_1)\times\mathcal H(G_2):(H_1,H_2)\text{ is good}\}
}{
|\mathcal H(G_1)|\,|\mathcal H(G_2)|
},
\]
with the convention $A_n(G_1,G_2)=0$ if either graph has no Hamilton cycle.  Thus \cref{thm:typical-hamiltonian-alignment} is the statement
\[
A_n(G_1,G_2)\xrightarrow{\mathbb P}1.
\]

\begin{proof}[Proof of \cref{thm:typical-hamiltonian-alignment}]
Use the marked planted law $\widehat{\mathcal Q}_{n,d}$ from \cref{subsec:ham-planted}, with unmarked graph marginal $\mathcal Q_{n,d}$. Conditional on the resulting graph $G$, the distinguished Hamilton cycle is uniform over $\mathcal H(G)$. Indeed, before conditioning, every pair $(H,R)$ with $H$ a Hamilton cycle and $R$ a simple $(d-2)$-regular graph has the same probability. After conditioning on edge-disjointness, a fixed $d$-regular graph $G$ and a Hamilton cycle $H\subseteq G$ determine exactly one residual graph,
\[
E(R)=E(G)\setminus E(H).
\]
Hence
\[
\Prob{H=H_0\mid G}=\frac1{|\mathcal H(G)|}
\]
for every $H_0\in\mathcal H(G)$. For two independent planted samples, the distinguished cycles are independent conditional on $(G_1,G_2)$.

Set
\[
B_n:=1-A_n.
\]
Under $\mathcal Q_{n,d}^{\otimes2}$, $B_n(G_1,G_2)$ is exactly the conditional probability that the two distinguished Hamilton cycles form a bad pair. The proof of \cref{thm:hamiltonian-rotation-overlap} showed that the distinguished pair is bad with probability $o(1)$. Therefore
\[
\mathbb E_{\mathcal Q_{n,d}^{\otimes2}}[B_n]=o(1).
\]
For every fixed $\delta>0$, Markov's inequality gives
\[
\mathcal Q_{n,d}^{\otimes2}(B_n>\delta)
\le
\frac{\mathbb E_{\mathcal Q_{n,d}^{\otimes2}}[B_n]}{\delta}
\longrightarrow0.
\]
The event $\{B_n>\delta\}$ depends only on the two unmarked graphs. Since
\[
\mathcal Q_{n,d}^{\otimes2}\approx\mathcal U_{n,d}^{\otimes2},
\]
contiguity yields
\[
\mathcal U_{n,d}^{\otimes2}(B_n>\delta)\longrightarrow0.
\]
As this holds for every $\delta>0$,
\[
B_n\xrightarrow{\mathbb P}0,
\qquad
A_n\xrightarrow{\mathbb P}1.
\]

For the equivalent sampling formulation, note that $0\le B_n\le1$. Hence
\[
\mathbb E_{\mathcal U_{n,d}^{\otimes2}}[B_n]
\le
\delta+\mathcal U_{n,d}^{\otimes2}(B_n>\delta).
\]
Letting first $n\to\infty$ and then $\delta\downarrow0$ shows that this expectation tends to zero. It is precisely the failure probability obtained by first sampling $G_1,G_2$ and then selecting one Hamilton cycle uniformly from each graph.
\end{proof}

Note that the theorem permits exceptional pairs of Hamilton cycles, but their proportion tends to zero in probability. It does not claim that every pair is successful, nor does it automatically apply to Hamilton cycles selected according to an arbitrary graph-dependent rule.

\subsection{Sharpness of the matching scale}
\label{app:ham-sharpness}

\begin{proof}[Proof of \cref{prop:cyclic-matching-overlap}]
The lower bound is \cref{lem:cyclic-matching-witness}. For the upper bound, fix $\varepsilon>0$ and set
\[
k=\left\lceil(1+\varepsilon)\frac{\log n}{\log\log n}\right\rceil.
\]
Condition on $J_1$. For any fixed shift $s$, every $k$-subset of $\rho_s(J_1)$ is a matching, and by \cref{lem:matching-containment} it is contained in $J_2$ with probability $p_k^{(a_2)}$. Therefore
\begin{align*}
\Prob{T_n\ge k\mid J_1}
&\le n\binom{a_1}{k}p_k^{(a_2)}\\
&\le\frac{n}{k!}\left(\frac{a_1}{n}\frac{2a_2}{n}\right)^k\exp(O(k^2/n))\\
&\le\frac{n\,2^{-k}}{k!}\exp(O(k^2/n)).
\end{align*}
By Stirling's formula,
\[
\log n-k\log k+O(k)=(-\varepsilon+o(1))\log n,
\]
so the last expression tends to zero. Together with \cref{lem:cyclic-matching-witness}, this proves the proposition uniformly over the allowed matching sizes.
\end{proof}

\begin{proof}[Proof of \cref{thm:cyclic-matching-unions}]
For the lower bound, apply \cref{prop:cyclic-matching-overlap} to $M_{1,1}$ and $M_{2,1}$.

For the upper bound, condition on $R_1$, which contains at most $q_1n/2$ distinct edges. Let $S$ be a prescribed collection of $k$ distinct edges. If $S\subseteq R_2$, each edge of $S$ belongs to at least one of the $q_2$ perfect matchings. There are at most $q_2^k$ ways of assigning the edges of $S$ to these matchings. For a fixed assignment, the probability is zero unless the edges assigned to each matching form a matching. Otherwise, independence and \cref{lem:matching-containment} give
\[
\Prob{S\subseteq R_2}
\le q_2^k n^{-k}\exp(O(k^2/n)).
\]
A union bound over shifts and $k$-subsets of $\rho_s(R_1)$ therefore gives
\[
\Prob{\max_s|\rho_s(R_1)\cap R_2|\ge k\mid R_1}
\le
\frac{n(q_1q_2/2)^k}{k!}\exp(O(k^2/n)).
\]
For
$
k=\left\lceil(1+\varepsilon)\frac{\log n}{\log\log n}\right\rceil,
$
the right-hand side is $n^{-\varepsilon+o(1)}$. Combined with the lower bound, this proves the theorem.
\end{proof}

\section{Technical estimates for the cycle-factor lower bound}\label{app:cycle-factor-technical}

This appendix collects the probabilistic and coupling estimates used in \cref{sec:lower}.  The deterministic star construction remains in the main text.  Throughout, $d$ is fixed.

\subsection{Local sparsity}\label{app:cf-local-sparsity}

\begin{lemma}[Local sparsity]\label{lem:cf-local-sparsity}
Let $d\ge3$ be fixed and suppose that
\[
k\le \frac{2}{a_d}\log n.
\]
Then there exists a constant $J=J_d$ such that,
with probability $1-o(1)$, a uniform simple
$d$-regular graph $G$ on $[n]$ satisfies
\[
e_G(S)\le |S|+J
\qquad
\text{for every }S\subseteq[n]\text{ with }|S|\le2k.
\]
\end{lemma}

\begin{proof}
We first work in the $d$-regular configuration model. By assumption, we may take $C=2/a_d$, so that $k\le C\log n$.

Fix a collection of $m=O(\log n)$ distinct edges. Each edge has two endpoints, and each endpoint has at most $d$ half-edges. Hence, there are at most $d^{2m}$ ways to assign the edge endpoints to configuration half-edges. For any compatible assignment, the probability that the resulting $m$ prescribed half-edge pairs occur in a uniform pairing is
\[
\prod_{j=0}^{m-1}\frac1{dn-(2j+1)}
\le
\left(\frac{2}{dn}\right)^m
\]
for all sufficiently large $n$.  Hence
\begin{equation}\label{eq:cf-prescribed-edge-bound}
\Pp(\text{all $m$ prescribed ordinary edges occur})
\le
\left(\frac{2d}{n}\right)^m.
\end{equation}

Fix $1\le v\le2k$.  If a $v$-vertex set spans at least $v+J$ distinct edges, then some prescribed collection of $v+J$ edges is present inside it.  By~\eqref{eq:cf-prescribed-edge-bound} and a union bound, the probability that such a set exists is at most
\[
\binom nv
\binom{\binom v2}{v+J}
\left(\frac{2d}{n}\right)^{v+J}.
\]
Using $\binom ab\le(ea/b)^b$ and $v+J\ge v$, this is bounded by
\[
(e^2d)^v\left(\frac{edv}{n}\right)^J.
\]
Summing over $1\le v\le2C\log n$ gives at most
\[
2C\log n\,(e^2d)^{2C\log n}
\left(\frac{2edC\log n}{n}\right)^J.
\]
Choosing the fixed integer $J=J_d$ sufficiently large makes this $o(1)$.

Thus, the local-sparsity event holds with probability $1-o(1)$ in the configuration model. Since the simplicity probability tends to a positive constant, the same high-probability statement holds after conditioning on simplicity.
\end{proof}

\subsection{Internal edges in a uniform pairing}\label{app:cf-pairing-concentration}

\begin{lemma}[Pairing concentration]\label{lem:cf-pairing-concentration}
Let $r\ge1$ be fixed and let $R$ be the configuration-model $r$-regular multigraph on $n$ vertices. If $U\subseteq[n]$ is fixed with $u:=|U|\ge\eta n$, where $\eta>0$ is fixed, then for all sufficiently large $n$,
\begin{equation}\label{eq:cf-pairing-concentration}
\Pp\left(
e_R(U)<\frac{ru^2}{4n}
\right)
\le
\exp\left(-\frac{r\eta^4}{256}n\right).
\end{equation}
\end{lemma}

\begin{proof}
There are $L=rn$ half-edges, of which $a=ru$ belong to the vertices of $U$. Let $X=e_R(U)$. Every unordered pair of half-edges in $U$ is matched with probability $1/(L-1)$, and therefore
\[
\Ee X
=\frac{\binom a2}{L-1}
=\frac{a(a-1)}{2(L-1)}
=\frac{ru^2}{2n}+O_r(1).
\]
Expose the uniform pairing one pair at a time, always revealing the half-edge matched to the first unmatched half-edge in a fixed ordering, and let $(Z_j)$ be the martingale of the final value of $X$. Two possible partner choices for the exposed half-edge can be coupled by switching at most two pairs. Such a switch changes the number of pairs lying entirely in $U$ by at most two, so
\[
|Z_j-Z_{j-1}|\le2.
\]
There are $L/2=rn/2$ exposure steps. For sufficiently large $n$, using $u\ge\eta n$,
\[
\Ee X-\frac{ru^2}{4n}
\ge
\frac{r\eta^2n}{8}.
\]
The Azuma--Hoeffding inequality (see \cite{mcdiarmid1989bounded}) therefore gives
\begin{align*}
\Pp\left(X<\frac{ru^2}{4n}\right)
&\le
\exp\left(
-\frac{(r\eta^2n/8)^2}{2(L/2)\,2^2}
\right)\\
&=
\exp\left(-\frac{r\eta^4}{256}n\right),
\end{align*}
which proves~\eqref{eq:cf-pairing-concentration}.
\end{proof}

\subsection{Deletion and repair}\label{app:cf-deletion-repair}

We now prove the size-reduction coupling used to remove the divisibility condition. The construction is easiest at the pairing level.

\begin{lemma}[Deletion-and-repair coupling]\label{lem:cf-deletion-repair}
Fix $d$.  Let $N<n$ be admissible and set $r=n-N=O_d(\log n)$.  There is a coupling of a uniform simple $d$-regular graph $G$ on $[n]$ and a uniform simple $d$-regular graph $H$ on $[N]$ such that, with probability $1-o(1)$,
\begin{equation}\label{eq:cf-coupling-edge-loss}
|E(H)\setminus E(G[[N]])|
\le\frac{dr}{2}.
\end{equation}
\end{lemma}

\begin{proof}
First, couple the unconditioned configuration models.  Sample a uniform pairing $P$ of the $dn$ half-edges.  Delete the $dr$ half-edges belonging to vertices $N+1,\ldots,n$, together with all pairs incident to them.  Let $A$ be the set of surviving half-edges whose matched half-edges were deleted, and write $a=|A|$.  Then $a$ is even and $a\le dr$. Conditional on the retained partial pairing, pair the elements of $A$ uniformly among themselves, and denote the resulting pairing on the $dN$ surviving half-edges by $P'$.

The distribution of $P'$ is invariant under every permutation of the $dN$ surviving distinguishable half-edges.  Indeed, the original pairing is invariant, and the deletion-and-repair operation is equivariant.  Since the permutation group acts transitively on perfect matchings of these half-edges, $P'$ is a uniform pairing.  Every pair of $P$ supported entirely on $[N]$ is retained, and exactly $a/2\le dr/2$ repair pairs are added.  Thus~\eqref{eq:cf-coupling-edge-loss} holds deterministically at the pairing level.

It remains to pass to simple graphs simultaneously. Let $S$ and $S'$ be the events that $P$ and $P'$ are simple.  We claim
\begin{equation}\label{eq:cf-simplicity-symmetric-difference}
\Pp(S\triangle S')=o(1).
\end{equation}
A loop or a parallel-edge defect of $P$ involving at least one deleted vertex has probability $O_d(r/n)=o(1)$ by the standard first-moment bounds for defects in the configuration model.  Every defect supported entirely on $[N]$ survives unchanged in the retained partial pairing.

New defects can only be created by the repair. Conditional on $a$, symmetry makes $A$ a uniform $a$-subset of the surviving half-edges. Hence, the expected number of pairs of half-edges in $A$ belonging to the same vertex is
\[
O_d\left(\frac{a^2}{N}\right)
=O_d\left(\frac{r^2}{n}\right)=o(1).
\]
With probability $1-o(1)$, therefore, all half-edges in $A$ belong to distinct vertices. Conditionally on this event, the repair creates no loops, and distinct repair edges cannot be parallel.

It remains to exclude parallel edges between the repair pairing and the retained pairing. Condition on $A$ and on the repair pairing, but not on the retained pairing. Given $A$, the retained pairing is uniform on the $dN-a$ remaining half-edges and is independent of the repair pairing.

Fix a repair edge joining distinct vertices $u$ and $v$. Since each vertex has at most $d$ half-edges in the retained pairing, the conditional probability that this pairing already contains an edge between $u$ and $v$ is at most
\[
\frac{d^2}{dN-a-1}=O_d(N^{-1}).
\]
There are $a/2\le dr/2$ repair edges, so a union bound gives
\[
\mathbb P(\text{a repair edge duplicates a retained edge})
=O_d\left(\frac{r}{N}\right)=o(1).
\]
Consequently, the repair creates no new defects with probability $1-o(1)$. Together with the bound on defects involving deleted vertices, this proves \eqref{eq:cf-simplicity-symmetric-difference}.

Both $\Pp(S)$ and $\Pp(S')$ tend to the same positive simplicity constant.  Since $\Pp(S\triangle S')=o(1)$, conditioning the joint construction on $S\cap S'$ changes either conditional simple marginal by total variation $o(1)$.  A maximal-coupling correction therefore yields exact marginals $G\sim\mathcal U_{n,d}$ and $H\sim\mathcal U_{N,d}$ while preserving~\eqref{eq:cf-coupling-edge-loss} outside an event of probability $o(1)$.
\end{proof}

\begin{proof}[Proof of \cref{lem:size-reduction-coupling}]
Take two independent copies of \cref{lem:cf-deletion-repair}.  Thus, for $i=1,2$, we couple $G_i\sim\mathcal U_{n,d}$ with $H_i\sim\mathcal U_{N,d}$ so that, with probability $1-o(1)$,
\[
|E(H_i)\setminus E(G_i[[N]])|
\le\frac{dr}{2}.
\]
The two couplings are independent, so $(G_1,G_2)$ and $(H_1,H_2)$ have the required product marginals.

Let $\phi\in\mathfrak S_N$ attain $M(H_1,H_2)$ and extend it to a permutation of $[n]$ by mapping the deleted vertices $N+1,\ldots,n$ among themselves.  Of the common edges counted by $M(H_1,H_2)$, at most $dr/2$ use an edge of $H_1$ absent from $G_1[[N]]$, and at most $dr/2$ use an edge of $H_2$ absent from $G_2[[N]]$.  Every other common edge remains common under the extended alignment.  Hence
\[
M(G_1,G_2)
\ge
M(H_1,H_2)-dr,
\]
which is~\eqref{eq:size-reduction-overlap}.
\end{proof}

\end{document}